\documentclass[reprint]{revtex4-2}

\usepackage{amsmath}
\usepackage{amsthm}
\usepackage{amssymb}
\usepackage{bm}
\usepackage{siunitx}
\usepackage{graphicx}
\usepackage{comment}
\usepackage[dvipsnames]{xcolor}
\usepackage{xspace}
\usepackage{algpseudocode}

\usepackage{algorithm}

\begingroup

\footnotetext{These authors contributed equally to this work.}
\endgroup
\begingroup

\footnotetext{Corresponding author: melvyn.tyloo@gmail.com}
\endgroup

\begin{document}

\author{Laurent Pagnier$^*$}
\affiliation{Department of Mathematics, The University of Arizona, Arizona, USA}

\author{Melvyn S.~Tyloo$^{*,\dagger}$}

\author{Akshita Jindal}
\author{Pragati Thakur}
\author{Kyle C.~A.~Wedgwood}

\affiliation{Living Systems Institute, University of Exeter, Exeter, UK }
\affiliation{ Department of Mathematics and Statistics, University of Exeter, Exeter, UK}
\title{A Closed-loop Framework to Discriminate Models Using Optimal Control}

\begin{abstract}
Predicting the response of an observed system to a known input is a fruitful first step to accurately control the system's dynamics.
Despite the recent advances in fully data-driven algorithms, the most interpretable way to reach this goal is through mechanistic mathematical modeling. 
Here, we leverage optimal control and propose a closed-loop iterative method to choose among a set of candidate models the one that most accurately predicts an observed system. 
We assume that one has control over an input of the observed system and access to measurements of its response. 
Our approach is to identify the input control that maximally discriminates the response of the candidate models, allowing us to determine which model is best by comparing such responses with the observed data.
We demonstrate our proposed framework in numerical simulations before applying it during an electrophysiology experiment, successfully discriminating between different models for photocurrents produced via opsin dynamics.
\end{abstract}

\maketitle

\section{Introduction}
Mathematical modeling is the primary tool used to predict and eventually control the future evolution of an observed system. 
However, finding the most adequate model is usually a challenging task, as there often exist multiple models with varying levels of complexity and hypothesized mechanisms for a given system.
One then has to fit these models to the available measurement data from the observed system, with the major issue that it may not be sufficient to distinguish the candidate models from one another~\cite{silk2014model}. 
This is typically the case when the data correspond to measurement of a macroscopic or output observable that depends on a multitude of microscopic degrees of freedom. 
For example, in biological systems, the underlying number of interaction degrees of freedom is usually large and the observability is often limited to e.g., firing rates, concentrations, overall output currents~\cite{chis2011structural,villaverde2019observability}\,.

Given some measurement data, it is tempting to apply equation discovery techniques~\cite{gerardos2025principled} to try to uncover the underlying dynamics at the origin of the data. This approaches suffers from the major drawback that the data might be limited to observation of the system within a subspace of the total phase space. Within such a subspace, several models might be equally likely to explain the data.
The output model from an equation discovery algorithm might therefore not be able to predict the future evolution of the system, especially if the experiment's conditions change.
Such limitation in the training data also prevents one from using fully data-driven machine learning algorithms to reproduce the data. Indeed, for such algorithms without any assumed model, predicting the behavior outside the phase space spanned by the training data set is very challenging \cite{brunton2022data}.

From another perceptive, model selection has been investigated using statistical tools, such as the Bayes factor, which compares the posterior probability of competing models. Other criteria like the Bayesian Information Criterion (BIC), based on the local geometry of the likelihood function around the maximum likelihood estimator (MLE) or the Akaike Information Criterion (AIK), leveraging the Kullback-Leibler (KL) divergence between the data and the candidate model have also been developed for model selection~\cite{ding2018model}.
These methods allow to compare how likely one model is over another to have generated some given data. However, as in the case of equation discovery discussed previously, if the data do not exhibit all the features of the actual underlying model, one might end up selecting an incorrect model.

Collecting data that enable one to choose one model over the other is of primal importance (i) to predict the behavior induced by inputs, which opens the way to further experimental designs and engineering applications; and (ii) to manage the cost of data acquisition by optimizing the information content of measurements. 
To achieve model discrimination and eventually select the most adequate model, one needs to design system inputs that generate observed outputs that allow candidate models to be distinguished. 
The so-called T-optimality criterion for designing static inputs for such a model discrimination task, which falls under the broader category of optimal experimental design, was first developed by Atkinson and Federov~\cite{atkinsonDesignExperimentsDiscriminating}.
This criterion, originally designed for linear models and later generalized to nonlinear models \cite{atkinson2010nonuniqueness,dette2009nonlinear}, designs inputs that maximize the error measured as the residual sum of squares between an alternative model and a `true model', so that the true model can be correctly selected as the most appropriate.
Other metrics, such as the KL divergence \cite{skandaOptimalExperimentalDesign2010} can also be used in place of the sum of squares to achieve a similar goal, and prior beliefs about which model is correct can be accounted for \cite{poncedeleon1991priorInfo}.
The approach may also be generalized to the case in which neither of the proposed models is true to provide a robust T-optimality criterion \cite{vajjah2012robustOptimality,detteBayesianToptimalDiscriminating2015}, with the latter reference adopting a Bayesian framework to account for model uncertainty.

For dynamical systems in which observations can be taken, and inputs applied, in a temporally continuous manner, as opposed to regression models considered in earlier studies, the question then becomes how to design inputs that bring the system to regions of the phase space where ambiguity between the models is ruled out. 
With this task in mind, we here propose a framework based on optimal control that designs system inputs which bring the system to a region of the phase space where candidate models can be distinguished.
We assume that one has access to a system that can be observed while applying different, arbitrary input control signals. 
Note that the choice of control signals to achieve this depends sensitively on the parametrization of the competing models, for which control inputs can also be designed~\cite{gadkar2005iterative,bandara2009optimal,kreutz2009systems}\,.
Such design of inputs for model parametrization can be achieved using the D-optimality criterion \cite{myung2009optimal}, which maximizes the determinant of the Fisher information matrix (FIM, which measures the sensitivity of the likelihood of a model under parameter variation), or the A- and E- optimality criteria, which maximize the trace or smallest eigenvalue of the inverse of the FIM, respectively (see Chapter 10 in \cite{atkinsonBookOptimumExpt}).

In a similar manner to the one we adopt for model discrimination, a continuous control signal can be used to steer system trajectories to regions of phase space for which properties of the FIM are maximized \cite{hooker2015control}.
These approaches can be also applied in a Bayesian sense using techniques from Bayesian experimental design, which is thoroughly reviewed in ~\cite{chalonerBayesianDesign}.
In general, the choice of inputs that achieve D- or A-optimality are not the same as those that achieve T-optimality, though attempts to achieve both goals simultaneously have been made, leading to the DT-optimality criterion\cite{atkinson2008dtoptimality,waterhouse2004dtoptimality}.
Moreover, when models are not identifiable (either structurally or practically) due to noise, parameter sloppiness, or partial observability, methods based on Fisher information may lead to inaccurate results \cite{liu2025optimal}.

Our general algorithm starts with a pool of candidate models, and sequentially fits their parameters and then, using optimal control, designs an input signal the maximizes the discrepancy between the models~\cite{murray2009optimization}. 
One then iteratively injects the obtained signal into the experiment and performs again model fitting and design of the next control signal.
For simple models, or ones that can be easily and robustly parametrized, sequential designs are not, in general, expected to be more efficient than non-sequential designs \cite{chalonerBayesianDesign}, however, this is not true when models are mis-specified \cite{myung2009optimal} or when the results of the previous recording may impact the design of the next set of control inputs, as in our case.
The process stops when one of the models is identified as the most adequate to predict the output of the experiment, or when a maximum number of iterations has been reached.
We note that Ref.~\cite{liu2025optimal} recently used optimal control using Pontraygin's Maximum Principle to design a signal that maximizes the difference between two models, under the assumption that models were identifiable and were fully parametrized.
While based on similar principles, our method goes one step further by re-injecting the designed signal into the observed system, and performing the model discrimination over multiple control realizations until one of the models can be selected over the other.

To showcase the algorithm and motivated by the increasing number of experiments using optogenetics that are performed in many research fields ~\cite{deisseroth2011optogenetics,emiliani2022optogenetics}\,, and in particular for modulation of neuronal circuits, we investigate opsin models. 
Such models have been introduced to describe the relationship between the output cellular photocurrent and the light input, including different levels of complexity and transition mechanisms.
While these models have a rather limited number of degrees of freedom, their observability and control are also very limited. 
Indeed, the intensity of the light shed on the cells is the only control input and the only model output that can be measured is the induced photocurrent. 
With a single input globally impacting the dynamics, a single output that is function of many variables, and many unobserved variables, finding the correct opsin model that best predicts the behavior of an optogenetics experiment is a highly complex task.
For this reason, we argue that the discrimination of such opsin models already represents a remarkable application of our method, for which our framework provides useful information both when using a numerically simulated system and during an actual laboratory experiment.
We also remark that while Bayesian frameworks may seem like the natural framework in which to study noisy biological systems, such approaches require considerable computational efforts to estimate likelihood functions (e.g., using Approximate Bayesian Computation), making them impractical for our system involving a closed-loop interaction between computation and data acquisition, even when using techniques that minimize the per iteration work required, for example, through importance sampling ~\cite{vanlierBayesianApproachTargeted2012a} or supervised classification \cite{hainy2022classification}. 

Approaches to `front load' the computational work through amortization also provides means to minimize per iteration work by optimizing the \textit{policy} for designing control signals, rather than optimizing the signals themselves \cite{fosterDeepAdaptiveDesign,ivanova2021idad,huang2025aline}.
Under this approach, design of the signals essentially amounts to sampling from the space of possible signals following a predetermined set of rules.
However, much like the optimality criteria defined earlier, this approach relies on having strong confidence that the models are well parametrized at the outset.
Since we do not assume that our models are initially well parametrized, we choose to adopt a closed-loop approach that updates parameters in an online fashion.
In turn, since computations are required to be fast and the system we are considering is not dominated by noise, we opt to use a deterministic framework, rather than a probabilistic one.


The manuscript is organized as follows. In Sec.~\ref{sec:meth}, we describe our method. More precisely, we introduce an algorithm for parameter fitting in Sec.~\ref{sec:param}. Then, in Sec.~\ref{sec:disc}, we introduce the optimal control framework used to design input signals. Finally, the iterative procedure allowing model discrimination is given in Sec.~\ref{sec:iter} and illustrated on opsin models using numerical simulation in Sec.~\ref{sec:num} and on an experimental setup in Sec.~{\ref{sec:exp}}. We conclude this article by discussing our method in Sec.~\ref{sec:discussion}.

\section{Methods}\label{sec:meth}
We aim to find, among a set of candidate models, the one that is most accurate at predicting the future behavior of our system of interest. We refer to the latter as the \textit{reference system} in the following. 
To achieve model discrimination, we perform a sequential pairwise comparison of the models, keeping at each step the best match. We elaborate on our choice for a pairwise comparison versus a global comparison below in Sec.~\ref{sec:disc}. 
We start with one pair of candidate models without any initial sorting. The first step is to fit their parameters to the measurement data generated by the reference system under some initial control input.

Fitting model parameters is typically a complicated task that often relies on optimizing a non-convex function. It strongly depends on the initial conditions of the search algorithms and the volume of the state space that the reference system explored in the data. 
In general, one should also be careful about the structural identifiability of the candidate models ~\cite{grewal1976identifiability,godfrey1985identifiability,villaverde2016structural}. 
Indeed, depending on how the parameters are defined, many different choice of parameter values might have the same input-output relation, making it impossible to uniquely recover the parameter values.
This is even more true if the input control only allows the reference system to explore a limited subset of the entire state space. 
Such structural non-identifiablity of the candidate models is not a major limitation of our algorithm. 
Indeed, models that yield the same output for multiple distinct sets of parameters belong to the same equivalence class. 
As such, because their output is the same, our algorithm will return the same control input, regardless of which set of parameters was chosen within the equivalence class. 
Here, we are not interested in finding unique sets of parameters, but rather in choosing the model that predicts the best the response of the reference system to an input signal.
The identifiability of the models is therefore not a requisite for our method, however, one will have to take that into account when checking the convergence of the algorithm. 
Generically, our algorithm requires the weaker condition that candidate models can be distinguished from each other.

\begin{figure*}[t!]
    \centering
    \includegraphics[width=0.95\linewidth]{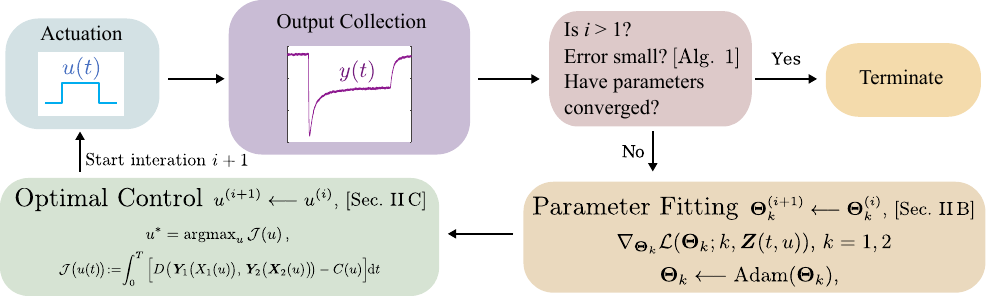}
    \caption{Schematic illustration of the closed-loop algorithm for model discrimination. For the first iteration, the initial control input $u(t)$ is arbitrarily chosen and produces the initial measurement dataset from the reference system. One then performs the parameter fitting of both candidate models ($k=1,2$). Using the fitted parameters, one maximizes the performance measure $\mathcal{J}$ over the control input $u(t)$\,. Finally, the obtained control input is used to actuate the reference system and the second iteration starts. This process continues until the models are successfully discriminated.}
    \label{fig:wf}
\end{figure*}

Following the fitting step, we solve an optimal control problem whose objective is to design an input control signal that maximizes the discrepancy between the two candidate model outputs.
The obtained control input is then used to obtain the next set of data from the reference system and evaluate how well each model predicted the output.
The algorithm then goes through another iteration until a termination criterion, which is discussed later, is reached.
One should note that iteratively updating the control input is essential to explore the state space and achieve a model discrimination that accounts for a large enough fraction of it. 
Moreover, the solution of the optimal control problem crucially depends on the fitted model parameters. 
Even if one of the candidate models is close to the actual one, poorly fitted parameters could lead to its non-selection.
Therefore, it is essential to fit the parameters at each iteration on the newly generated data. 
Heuristically, one expects this method to converge if the level of complexity of a least one of the proposed models matches the complexity of the observed system.
The overall workflow of the algorithm is summarized in Fig. \ref{fig:wf}.

\subsection{Model Formulation}
Let us give a more formal mathematical description of the dynamical systems we consider. 
We focus on specific systems where one has the control of an input signal $u(t)$\,, which is a scalar function in the following, but may be multivariate in general, and measures an output $\bm Y$ that is a function of the model variables and parameters.
Therefore, we consider models of the form
\begin{subequations}
\begin{align}
    \dot{\bm X} &= \bm F(\bm X,\, u(t);\, {\bm \Theta})\,,\label{eq:mod1}  \\
    \bm Y &= \bm H(\bm X;\, {\bm \Theta})\,,\label{eq:inf_y}
\end{align}
\end{subequations}
where ${\bm F}$ describes the evolution of state ${\bm X} \in \mathbb{R}^N$ subject to a parameter set ${\bm \Theta} \in \mathbb{R}^p$ and scalar, time-dependent control input $u(t)$; and ${\bm H}$ is the mapping from ${\bm X}$ to an observable ${\bm Y}\in \mathbb{R}^{M}$, which may also depend on ${\bm \Theta}$.
We assume that we have access to data ${\bm Z}(t;u) \in \mathbb{R}^{M}$ that are generated under forcing of a reference system with the control input $u(t)$.

In Sec. \ref{sec:num}\,, we numerically simulate the dynamics and collect the data, such that the actual model for the reference system is known. By contrast, in Sec. \ref{sec:exp}\,, the reference system is a laboratory experiment, such that there is no a priori known model for the data.

\subsection{Parameter Fitting}\label{sec:param}
We now introduce the learning scheme that we use to update the model parameters at each iteration. 
To fit the parameters ${\bm \Theta}_k$ for each of our models, we define the following loss function, 
\begin{align}\label{eq:loss}
    \mathcal{L}(\bm \Theta_k;k)  = \int_0^T D\Big({\bm Z}(t),{\bm Y}\big(\bm X_k(t)\big) \Big)\,{\rm d}t\,,
\end{align}
where $D(\cdot,\cdot)\geq 0$ is a distance function and $T$ is the duration over which the data are collected. Many options are available for $D$.
Here, we choose the Euclidean distance, i.e., $D(\bm x, \bm y) = \| \bm x - \bm y \|_2^2$\,, and optimize over the parameters such that the time integral in minimized.
For small models, the gradient of loss function~(\ref{eq:loss}) is straightforwardly obtained by applying automated differentiation over the numerical integration. Larger models require the use of a more dedicated method, such as the adjoint method. 
Parameters are then updated using this obtained gradient according to a gradient descent method, e.g., Adam~\cite{kingma2014adam}.
The number of epochs $N_{\rm grad}$ is one of Adam's hyperparameters that has to be selected carefully to allow efficient parameter updates but prevent overfitting.

As a final note, measurement data from the reference experiment is typically recorded at a finite sampling rate and is hence intrinsically discrete in time, i.e., ${\bm Z}(t; u):=\{{\bm Z}(t_0),{\bm Z}(t_1), \ldots\}$.
We therefore assume that this sampling rate is high enough to accurately capture the continuous system dynamics, e.g., via polynomial interpolation.

\subsection{Model Discrimination}\label{sec:disc}
Our goal is to find, among the set of candidate models, the one that best predicts the observed data ${\bm Z}(t;u)$\,, which depend on the input control signal $u(t)$\,. 
One way to achieve that is to find the input signal that maximizes simultaneously the prediction of all the candidate models. 
Another is to proceed with a binary tournament where pairs of candidates models are compared and only the best one is kept at each round. 
For this tournament, each model could be assigned a pre-ranking (similar to a seed in knockout tournament) based on criteria such as i) evaluation of $\mathcal{L}({\bf \Theta}_k,k)$ to some initial control profile, e.g., the box input shown in Fig.~\ref{fig:wf}, ii) model complexity, and iii) expert opinion on model plausibility.
The rankings of each round's winners would then be updated using the terminal values of $\mathcal{L}$ for that round, promoting the best models to proceed to the final round.
The approach is likely to be suitable in situations where the number of candidate models is low and the pairwise comparison is computationally faster than the global one.
Moreover, adopting a binary tournament is likely to be optimal when the wet lab experiments themselves can be parallelized, for example, using high throughput, automated systems \cite{seibertz_modern_2022}.
When the number of candidate models is large, the option based on simultaneous comparisons is likely to be more efficient, and we include an example of how to perform this in Sec.~\ref{SimModDisc} of the Supplementary Materials (SM).
Because time is critical in the experiment that we consider below and the number of candidate models is limited, we focus here on the binary tournament setting.

We aim to discriminate between two candidate mathematical models that could potentially reproduce the observed data ${\bm Z}(t;u)$\,, which depend on the input control signal $u(t)$\,.
Each model has its variables, ${\bm X}_k$\,, parameters ${\bm \Theta}_k$ and observables ${\bm Y}_k$, $k=1,2$. Note that state and parameter spaces are not necessarily the same, but we require the observable space to be the same for both models to allow for a meaningful comparison.
Similarly, the control input signal $u(t)$ should be common between the candidate models and the reference system. 
We introduce a subscript notation to denote the dynamics of the two candidate models as
\begin{subequations}
\begin{align}
    \dot{\bm X}_k &= \bm F_k\big(\bm X_k,\, u(t)\,;\bm \Theta_k\big)\,,\label{eq:x1}\\
    \bm Y_k &= \bm H_k\big(\bm X_k;\, \bm \Theta_k\big)\,,\label{eq:y1} \qquad k=1,2.
\end{align}
\end{subequations}

We seek a new control input signal $u(t)$ that provides the best possible discrimination between the two models, i.e., the largest difference between the two observables ${\bm Y}_1$ and ${\bm Y}_2$.
We achieve this goal using optimal control. First, we define a performance measure whose maximization provides a control input that generates a large discrepancy between the candidate models. The performance measure is given by,
\begin{align}\label{cost}
u^* &= {\rm argmax}_{u}\,\mathcal{J}(u)\,,\\
    \mathcal{J}\big(u(t)\big) &:= \int_0^T \Big[D\Big(\bm Y_1\big(\bm X_1(u)\big),\, \bm Y_2\big(\bm X_2(u)\big)\Big)- C(u)\Big]\,{\rm d}t\,, \nonumber
\end{align}
where $D$ is the same distance function as for the loss function Eq.~(\ref{eq:loss}); $C(u)$ is a cost for the control input that aims to avoid unrealistic signals, e.g., with a range of variation that is too large. This optimal control input can be found by deriving the adjoint system -- with  ${\bm \lambda_1}$\,, ${\bm \lambda_2}$ Lagrange multipliers to ensure that the dynamics of each model is satisfied -- and then by solving the associated two-point boundary value problem. 
Here, we choose a different approach that relies on solving an optimization problem under constraints where the degrees of freedom of the models are treated as optimization variables.

In the scheme description, we used a continuous point of view as it gives the lightest and most readable formulation. However, for the implementation, it is advantageous to discretize the system. Assuming a small, fixed timestep, $\Delta t$, the trapezoidal rule gives
\begin{align}
&\int_0^T f(t)\, {\rm d}t \approx \Big[(f(t_0)+f(t_N))/2 + \sum_{m=1}^{N-1}f(t_m)\Big]\Delta t\,,\nonumber\\
&\text{ with } t_m=\frac{T}{N\Delta t} \cdot m\,, 
\end{align}
for any integrable function $f$. Representing the evolution of the function $f$ as a vector $\bm f:=\left(f_0, \dots, f_N)\right)^\top$, where $f_m = f(t_m)$ and defining the operator
\begin{equation}
\mathcal{I}(\bm f):=\Big[(f_0+f_N)/2 + \sum_{m=1}^{N-1}f_m\Big]\Delta t\,,
\end{equation}
the optimization problem \eqref{cost} at iteration $i > 1$ can be rewritten and extended as
\begin{subequations}
\begin{align}\label{opt}
    u^{(i)} &= {\rm argmax}_{\{\bm u\}}\; \mathcal{I}\Big(\big\|\bm Y_{1} - \bm Y_{2}\big\|^2\Big)    
     \nonumber\\
    &- c_1\,\mathcal{I}\big(\bm u^2\big) - c_2\sum_{j=1}^{J}\mathcal{I}\big(\bm u\circ \bm u^{(i-j)}\big) \,, \\
&\textrm{s.t.} \nonumber\\
&\bm X_{k,m+1} = \bm X_{k,m} + \Delta t\,\bm F_k(\bm X_{k,m}, u_m)\,,\; \quad k=1,2\,, \label{eq:c_x}\\
&\bm Y_{k,m+1} =\bm H_k(\bm X_{k,m})\,, \label{eq:c_y}\\
&u_{\min} \le \mathcal{I}(\bm u^2)\,,\label{eq:min}\\
&\quad 0 \leq u_m \leq u_{\max}, \quad\, \;\quad m=1,\ldots,N\,,\label{eq:bound}
\end{align}
\end{subequations}
where $\circ$ is the (Hadamard) entry-wise product and $\bm u^2=\bm u\circ \bm u$,  and $c_1,c_2, u_{\rm min}, u_{\rm max}>0$ are hyperparameters and $N = T/\Delta t$\,. The last term in Eq.~\eqref{opt} is added to promote variety in the control profiles; without it the scheme may tend to become stuck in a series of similar control profiles. Thus, in addition to an intrinsic cost for the amplitude of the control, we add a term that rewards discrepancy with the $J > 1$ previously used control input signals, $\{\bm u^{(i-1)},\ldots,\bm u^{(i-J)}\}$. Namely,
\begin{align}
\int_0^T\big(u - u^{(i-j)}\big)^2 {\rm d}t &\nonumber\\
\approx  \mathcal{I}(\bm u\circ \bm u) &- 2\,\mathcal{I}\big(\bm u\circ \bm u^{(i-j)}\big) + {\rm const.}
\end{align}

The desired balance between intrinsic and variety-promoting costs is obtained by selecting $c_1$ and $c_2$ carefully. This set of control profile is in practice a fixed size memory, where $J$ must remain small for the optimization to be solved efficiently.  Here, we set $J=5$. Constraints (\ref{eq:c_x} - \ref{eq:c_y}) enforce the models' dynamics using forward Euler scheme. Constraint \eqref{eq:min} is provided to encourage exploration of the system response under forcing, which typically helps the algorithm to convergence to a solution more rapidly and avoid local maxima where the control input is vanishing at all $t$. Finally, constraint \eqref{eq:bound} ensures that control is always bounded. This problem is solved using state-of-art interior point method Ipopt~\cite{wachter2006implementation} within the optimization environment JuMP in Julia~\cite{Lubin2023}. The resulting optimal control, $\bm u^{(i)}\leftarrow \bm u$, can then be used as a signal to generate a next dataset ${\bm Z}(t;\bm u^{(i)})$.


\subsection{Iterative Discrimination and Parameter Estimation}\label{sec:iter}

As previously mentioned, our method iteratively discriminates models and updates their parameters.  
At each iteration, the algorithm proceeds by producing the next sequence of control inputs $\bm u^{(i)}=\{ u_m^{(i)}\}_{m=1,\ldots, N}$.
Then, the system is actuated with this optimal control input and its trajectory ${\bm Z}(t;\bm u^{(i)})$ is recorded. 
Using the method described in Sec.~\ref{sec:param}, the model parameters $\bm \Theta_1$ and $\bm \Theta_2$ are updated. 
These steps are repeated until termination. 

As illustrated in Fig.~\ref{fig:wf}, the process runs as long as the error in the prediction of the models is too large and the parameters have not converged. 
What is meant by this is that the algorithm needs some time (i.e., a few iterations) to update the models and give them the best chances to capture the system dynamics before the final model discrimination is made. 
Hence, we introduce a termination criterion that account for both the model fitness and the parameter updates, see Alg.~\ref{alg:stop}, with a goal to stop when we are confident that one of the models is significantly better than the other, that this model predicts the system response sufficiently well, and that this is not due to overfitting, the latter of which we assess by examining how much parameter change under variation of the control signals.
Thus to stop conclusively, the parameters of at least one of the models must have converged and the error between the prediction of that model and the data must be below a prescribed tolerance. 
We here note that convergence of parameters does not guarantee that the algorithm has found the globally optimal values due to prospective issues with structural and/or practical identifiability~\cite{wieland_structural_2021}.
\begin{figure*}[t!]
    \centering
    \includegraphics[width=0.9\textwidth]{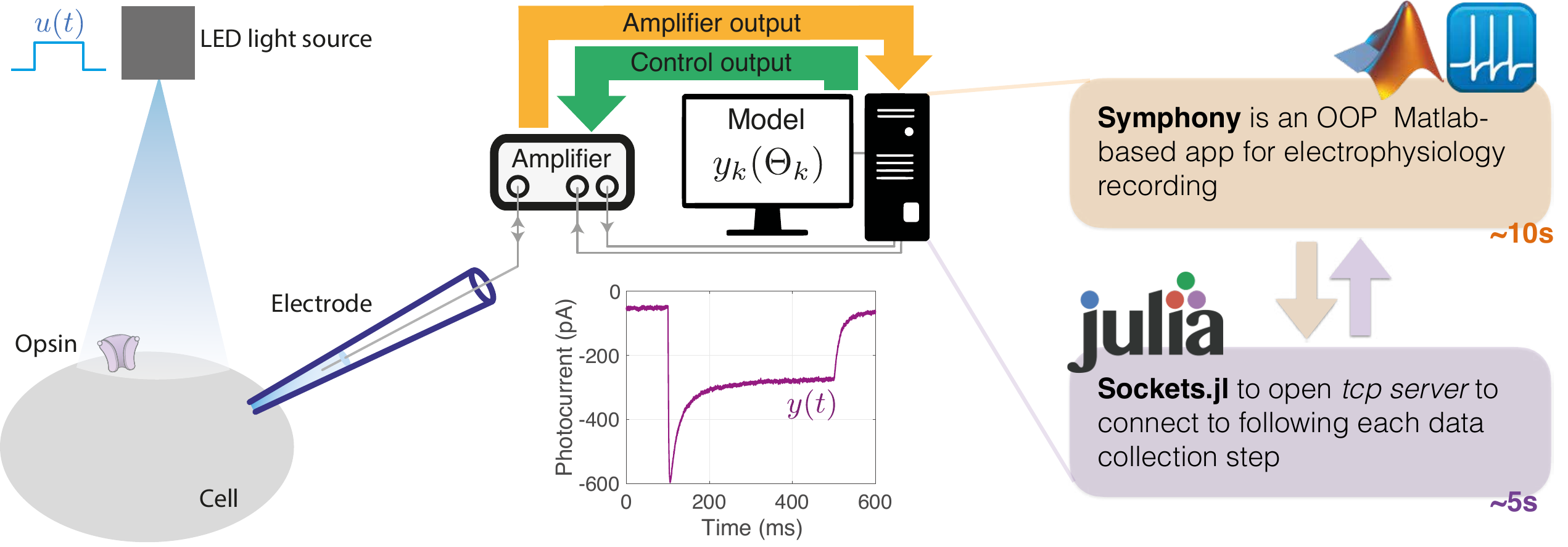}
    \caption{Schematic diagram of induction and measurement of photocurrents under voltage clamp. The left portion of the figure illustrates the experimental setting including hardware and biological components, the middle portion illustrates the flow of information between the hardware and software during a closed-loop experiment, and the right portion summarizes the software used.
    The graph in the lower-middle portion is an example of a photocurrent, which is our observable, under the application of the box input to control the LED intensity, indicated in the top-left.}
    \label{fig:patch}
\end{figure*}
\begin{algorithm}[H]

\begin{algorithmic}[1]
\Require Max steps $i_{\max}$, model fitting threshold $\mathcal{L}_{\max}$ and parameter increment threshold $\Delta_{\max}$
\If {$i > i_{\max}$ }
    \State stop as \textit{inconclusive};
\ElsIf {$\mathcal{L}_k < \mathcal{L}_{\max}\,,\;  k=1 \textbf{ xor } k=2 $ }
    \State $k^*\longleftarrow {\rm argmin}_k \mathcal{L}_k$;
    $\Delta {\bm\Theta}_{k^*} = {\bm\Theta}_{k^*}^{(i)} - {\bm\Theta}_{k^*}^{(i-1)}$
    \If { $\max_j|\Delta \Theta_{k^*,j}| < \Delta_{\max} $ }
        \State stop as \textit{conclusive} with $k^*$ as fittest;
    \Else
        \State continue discriminating the models;
    \EndIf
\ElsIf {$\mathcal{L}_k < \mathcal{L}_{\max}\,,\;  k=1 \textbf{ and } k=2 $ }
    \If { $\max_j|\Delta \Theta_{k^*,j}| < \Delta_{\max} $ }
        \State stop as \textit{conclusive} with the
         lowest complexity 
        \State model;
    \Else
        \State continue discriminating the models;
    \EndIf
\Else
    \State continue discriminating the models;
\EndIf
\end{algorithmic}
\caption{Termination Criterion}\label{alg:stop}

\end{algorithm}

Finally, we must handle the case where after the whole procedure the two models reproduce similarly the system dynamics. In this case, we apply Occam's Razor and select the model with the fewest parameters.

\section{Numerical Experiments}\label{sec:num}
To test and illustrate our method, we focus on models describing electrical currents generated in cells when light-sensitive proteins called opsins are activated by light. 
This type of light induced current has proven to be useful in a variety of electrophysiology applications~\cite{deisseroth2011optogenetics,emiliani2022optogenetics}. 
After introducing the opsin dynamics below, we first investigate numerically the method by simulating the reference system. 
We then apply the method in a closed-loop laboratory experiment.

\subsection{Model Description}\label{sec:opsin}

We now apply our scheme, displayed in Fig.~\ref{fig:wf}, to the mathematical modeling of photocurrents induced via light stimulation of channelrhodopsin.
Briefly, opsins are proteins that change conformation under illumination under specific light wavelengths.
In the case of channelrhodopsin (hereon referred to as opsin), this protein is a transmembrane ion channel that transitions from a non-conducting to a conducting conformation under stimulation of light of wavelength $\sim470$ nm.
This transition gives rise to a photocurrent that can be directly measured, for example, using the patch (voltage) clamp technique, depicted in Fig.~\ref{fig:patch}.
Here, an electrode is used to set the voltage across a cell membrane, in which many opsins are typically embedded, at a command value whilst currents across the membrane are measured under illumination of the cell at different intensities.

\begin{figure*}[t!]
    \centering
    \includegraphics[width=0.8\textwidth]{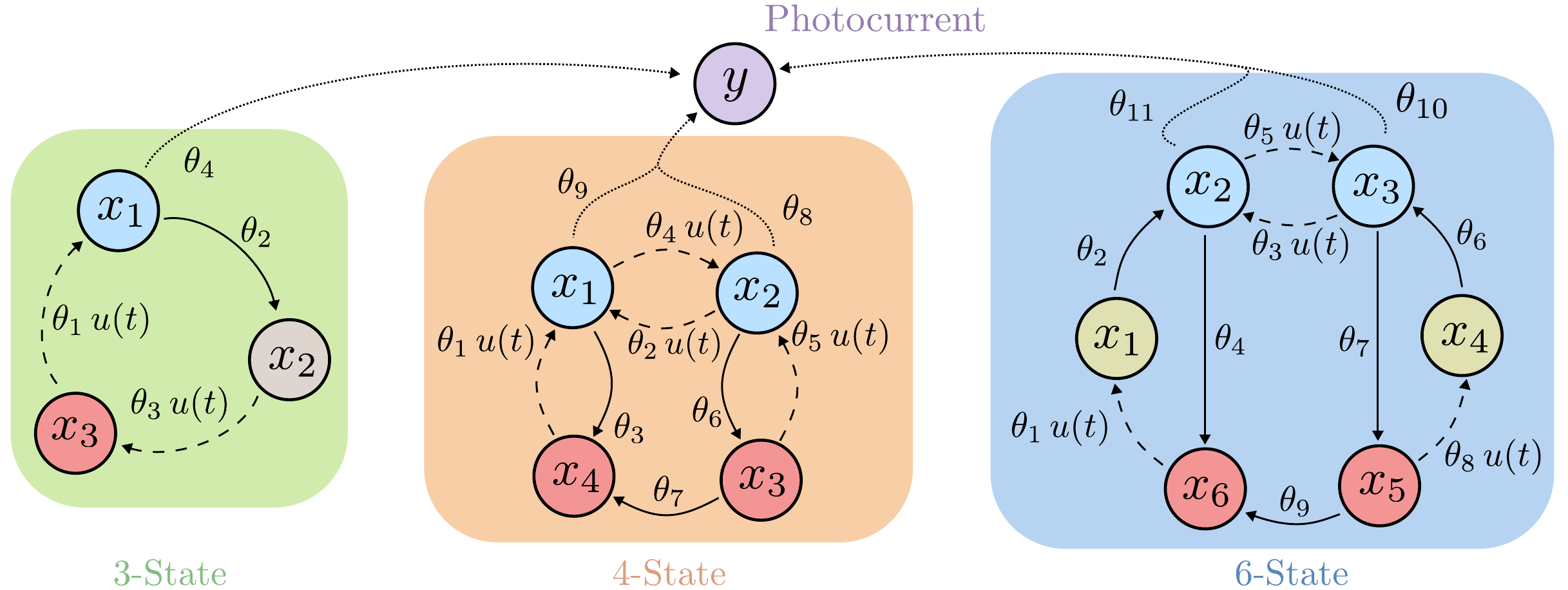}
    \caption{Depiction of 3-, 4-, and 6-state models. Blue, red, green and gray discs represent open, closed, intermediary and desensitized states respectively. Arrows show transitions and how they are parameterized and depend on control. In particular, dashed arrows indicate light-sensitive transitions; solid arrows represent light-insensitive transitions.}
    \label{fig:opsin}
\end{figure*}

Transitions between protein conformations are often described using simple $n$-state Markov models,  in which variables represent the probability that a protein is in a given state (as summarized in Fig.~\ref{fig:opsin}) and model terms represent the per unit time transition probabilities between states. Here, we use the `rate' version of such models, which assume that the number of proteins in the system is infinite.
Under this assumption, variables represent the proportion of channels in a given state and model terms represent transition rates between states.
Three distinct Markov models -- a 3-state, a 4-state, and a 6-state -- were put forward in \cite{evans2016pyrho}, as summarized in Fig.~\ref{fig:opsin}, to describe the dynamics of the channelrhodopsin we consider and so we select these as our prospective models.

In the simplest description -- the 3-state model -- which was first introduced in~\cite{nikolicPhotocyclesChannelrhodopsin22009}, each opsin channel can be in one of three states: 1) an open state in which ions can freely pass through and hence current can flow; 2) a closed state in which ions cannot pass; and 3) a desensitized state in which the channel is closed and cannot be opened.
Changes to the proportion of opsins in each of these states is then modeled as
\begin{subequations}
\begin{align}
\dot x_1 &= \theta_1  u(t)\,x_3 - \theta_2  x_1\,,\label{3s1} \\
\dot x_2 &= \theta_2  x_1 - \theta_3 u(t)\, x_2\,,\\
x_3 &= 1 - x_1 - x_2\,, \label{3sConstraint} \\
y &= \theta_4 x_1\,, \label{3s2}
\end{align}
\end{subequations}
where $x_1$ is the proportion of opsins in the open state, $x_2$ is the proportion of opsins in the desensitized state, and $x_3$ is the proportion of opsins in the closed state.
Here, we will always assume that the system was not excited prior to each experiment or, equivalently, that it had time to relax to its steady state, here $[0,0,1]^\top$.

Note that since the $x_i$ variables represent the proportion of opsins in a given state, conservation of mass implies that they must sum to 1, as enforced by \eqref{3sConstraint} in the above.
Hence, an $n$-state model has $n-1$ degrees of freedom.
The $\theta$ parameters govern the transition rates between states, a subset of which are scaled by the control inputs.
The control input, which represents the intensity of light~\footnote{Stricto sensu, the control is the voltage supplying the light which is the quantity that is measure and set, but here we will assume that there is a one-to-one correspondence between the two.} to which the opsins are exposed, affects a subset of these transition rates, and hence biases transitions to particular states. For sake of readability, we have relegated the description of the 4-state and 6-state models to the SM, Sec.~\ref{opsins}, but we note here that the 6-state model has an additional type of state, the \textit{intermediary state}, to account for second order dynamics in the transition between closed and open states (as introduced in~\cite{grossmanSpatialPatternLight2013}).
In theory, one would like to keep the light intensity not too large, as it may lead to physiological changes in the cells, or may even damage the cell membrane due to heating.
In practice, one can easily limit the maximum intensity such that even when the upper limit is hit, the physiological impact on the cells is minimal.
The experiments in Sec.~\ref{sec:exp} use an analog control voltage in the range 0V--10V to control an LED system.
As such, we set $\left( u_{\min}, u_{\max} \right) = (0,10)$ in what follows.
We further assume that we have fine temporal control of this analog voltage, unlike the bang-bang control used in \cite{liu2025optimal} since we see  distinct and reproducible responses as we vary the profile of the voltage input.

We now consider different cases using synthetic data, generated by one of the opsin models as the reference system.

\subsection{Validation in the Simplest Cases}

\begin{figure}[h!]
    \centering
    \includegraphics[width=\columnwidth]{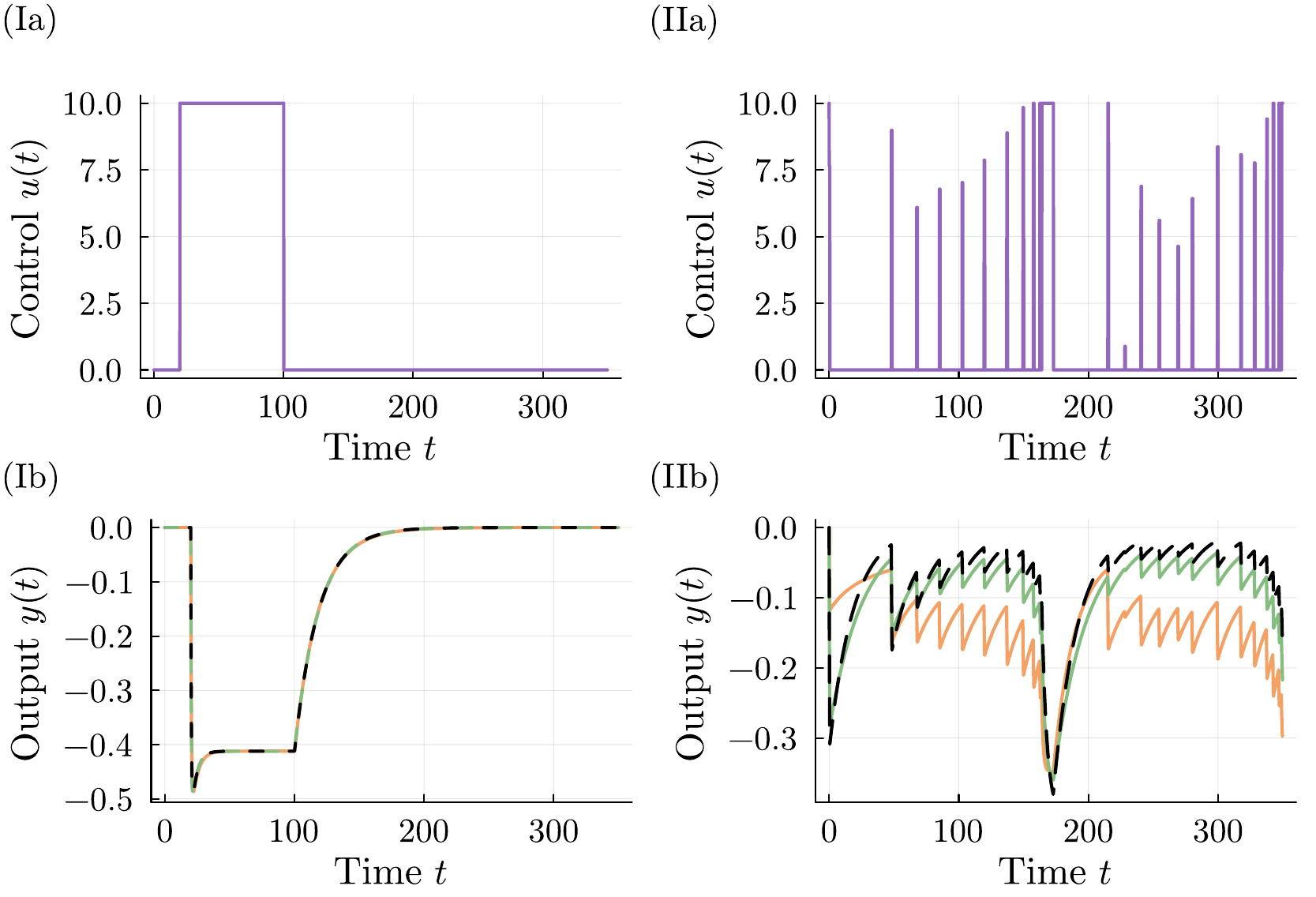}
    \caption{Box control input and parameter fitting (I), i.e. the standard procedure, compared to our procedure (II). (a) Control profiles. (b) Corresponding responses of the 3-state (green) and 4-state (orange) models reference model (dashed black), which in this case is the 3-state model with tuned parameters.
}
    \label{fig:case_2}
\end{figure}

For now, we will assume that the reference system has the same structure as one of the two candidate models that we aim to discriminate. In a first scenario, the candidate models are the 3-state and 4-state models while the reference system is simulated using the 3-state model. 
We show that, using an input control signal that is not diverse enough prevents one from being able to discriminate between different models. 
In particular, Fig.~\ref{fig:case_2} (Ia) and (Ib) show that using the standard procedure to characterize photocurrents, namely, to excite the system with a box input signal, is not suitable for model discrimination as the 4-state model is capable of overfitting to this single trajectory and so the predictions of the two models overlap.

On the other hand, Fig.~\ref{fig:case_2} (IIa) and (IIb) show that, when the input control has been obtained though our optimal control strategy, the 3-state model is more accurate in reproducing the reference system's trajectory. 
Indeed, Fig.~\ref{fig:case_2} (IIa) and (IIb) show, respectively, the control and trajectory at the end of our procedure. 
The 3-state model now outperforms the 4-state one in predicting the response of the reference. 
Note that the parameters in both models have been updated by the iterative parameter fitting process with more challenging control  and now the 4-state model is not able to (over)fit the system trajectory anymore.

\begin{figure}[h!]
    \centering
    \includegraphics[width=\columnwidth]{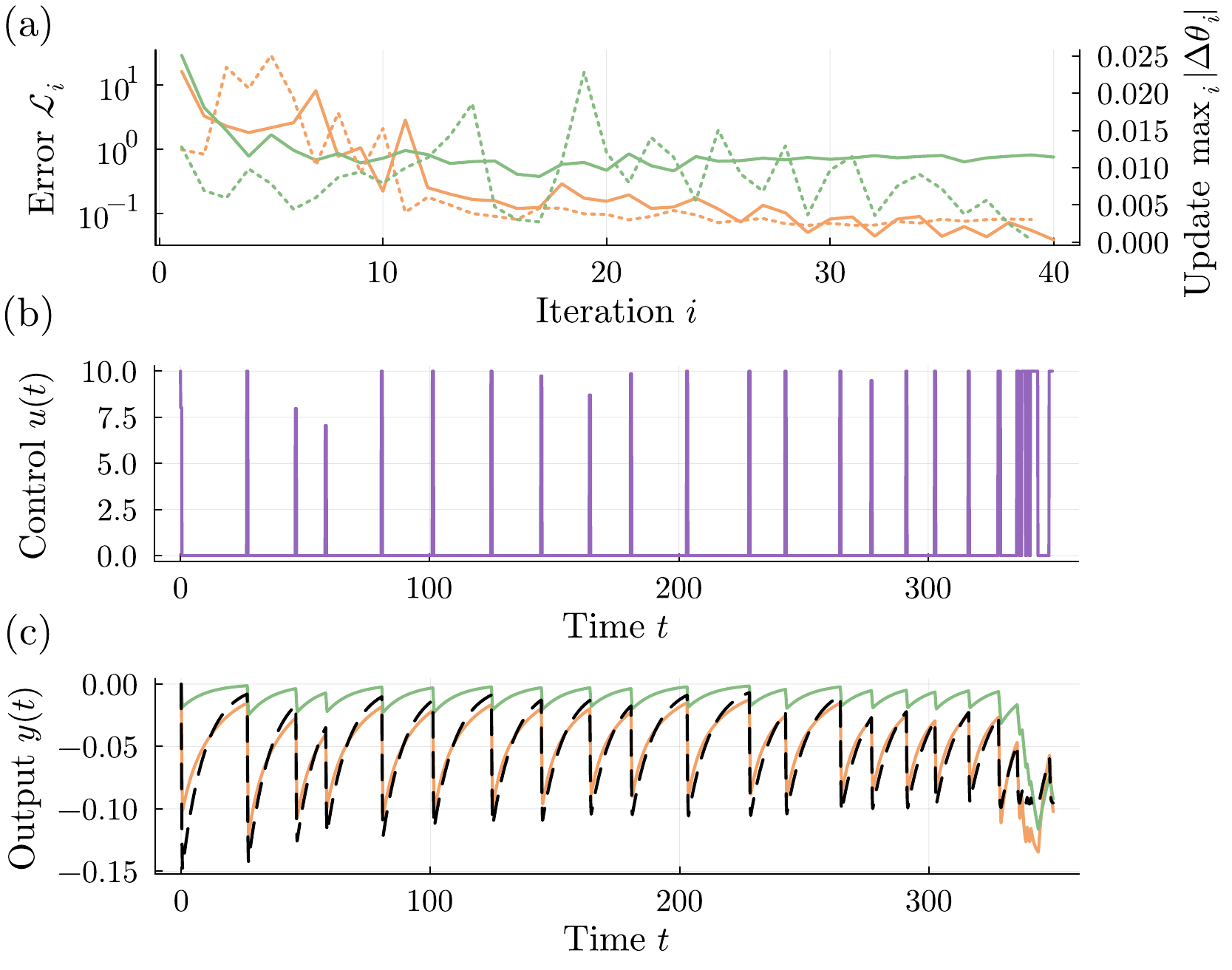}
    \caption{(a) Evolution of mismatch with the reference trajectory (solid) and of maximal parameter increment (dotted) for the 3-state (green) and 4-state (orange) models. (b) Optimal control at the end of the procedure. (c) Corresponding responses of both models compared with the reference model (dashed black), which in this case is the 4-state model with tuned parameters.}
    \label{fig:case_1}
\end{figure}

In a second scenario, the roles are switched, i.e., the reference system is a 4-state model and the two candidate models are the 3-state and 4-state. In Fig.~\ref{fig:case_1} (a)\,, we monitor the fitness of the models and observe that, at the start, both models (the parameters of which were selected randomly) do not capture the system's dynamics at all and large parameter updates are made. 
After a few more iterations, the 4-state model starts to perform better than the 3-state with a lower error, but it is to early to declare it as the winner as the mismatch with the reference system is still important and parameter updates are still large. 
Spikes in Fig.~\ref{fig:case_1} (a) happen when the optimal control changes abruptly and drives the system to configurations that have not been explored yet and for which the models are a poor fit. 
Finally, the 4-state model outperforms its competitor and its parameters converge. This behavior corroborates the termination criterion defined in Alg.~\ref{alg:stop}. Fig.~\ref{fig:case_1} (b) and (c) show the final optimal control and the corresponding outputs of both models. 
This provides a visual confirmation that the 4-state model is indeed the most accurate model to predict the response of the reference system to an input.

We conclude that our method correctly discriminates between the models in an ideal setup. We next move to the more challenging situation where the observed system is not exactly represented by any of the proposed models.

\subsection{Validation in More Complex Cases}\label{sec:complex}

We next apply our model discrimination method between every combination of the candidate models and reference system. 
The resulting outcomes are summarized in Tab.~\ref{tab:comp}\,, which gives the prediction errors obtained for each model discrimination. 
For each choice of reference system, one observes that the candidate model corresponding to the reference system always outperforms the other ones. 
\begin{table}[h!]
\includegraphics[width=\columnwidth]{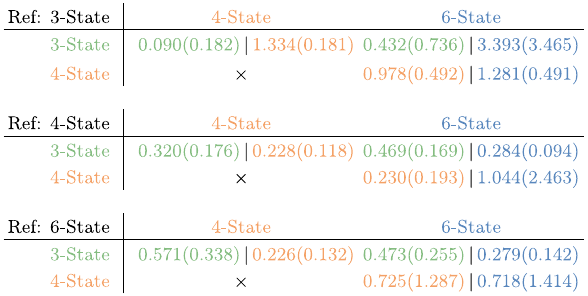}
\caption{Each model is used as a reference system and each pair of models is discriminated. Entries display the average discrepancies, with standard deviations in parentheses, at the last of 40 algorithm iterations. Averaging is performed over 10 samples with random initialization of the models' parameters. Colors indicate the model with which values are associated.}\label{tab:comp}
\end{table}
Our method therefore correctly identifies the most adequate model to predict the output of the reference system. 
Interestingly, when the reference system has not been simulated with one of the candidate models, sometimes one still obtains rather small prediction error. 
In a situation where the actual model for the reference system is not known, this could lead to the selection of the most accurate model to predict the response, even if the latter is not exactly the one of the reference. 
When the error remains too large, then one cannot draw any conclusion.

So far, we have tested our closed-loop algorithm on synthetic, noiseless data. 
However, most measurements of real world processes include some level of uncertainty. 
Although computational speed requirements lead us to adopt a deterministic framework for the fitting and control design steps, it is pertinent to address whether this approach is robust to noise in the reference system, which may appear as measurement noise in \eqref{eq:inf_y}, process noise in \eqref{eq:mod1}, or a combination of the two. We examine these two contributions separately. 

\begin{figure}[h!]
    \centering
    \includegraphics[width=\columnwidth]{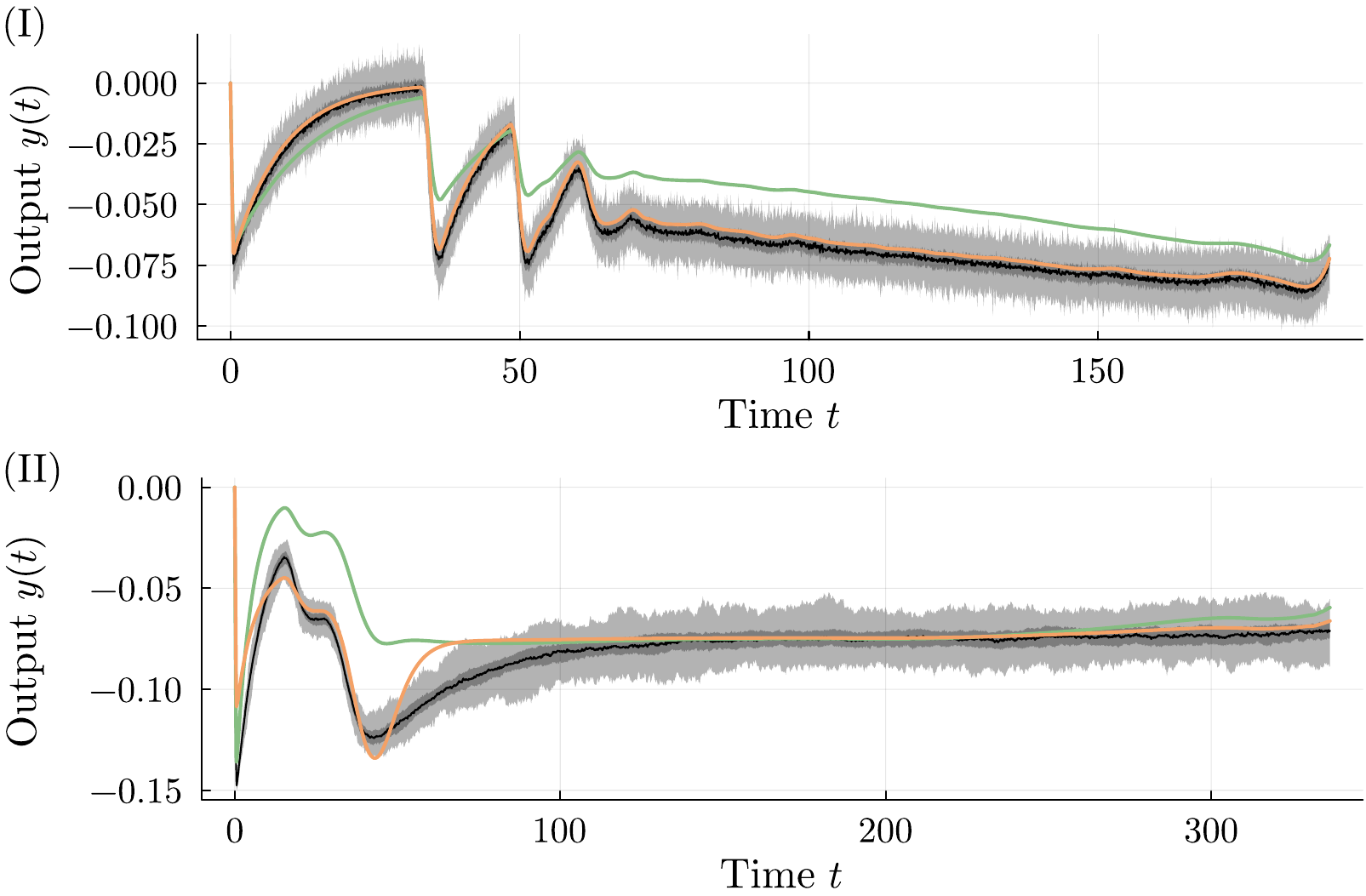}
    \caption{Two experiments with measurement noise (I) and stochastic dynamics (II). Reference model is the 4-State model in both cases. Trajectories of the 3-State (green) and 4-State (orange) models are compared to the reference (grey). As the references are stochastic models, their evolutions are displayed as quantiles: 1-99\% (light grey), 25-75\% (darker grey) and average (black line).}
    \label{fig:noise}
\end{figure}

We start with measurement noise. In this case, the trajectory $\bm Z(t,\bm u):= \{\bm Z(t_0),\ldots\}$ is obtained as noisy measurements,
\begin{equation}
\bm Z(t_i) = \bm H(\bm X^*(t_i),\bm \Theta^*) + \bm \xi_i\,,
\end{equation}
where $\bm \xi_i \sim \mathcal{N}(\bm 0, \bm \Sigma)$ and $\bm \Sigma$ is the covariance matrix that characterizes noise.  Fig.~\ref{fig:noise} (I) shows that our method is still able to discriminate the models accurately when the fitting target is noisy observation. The effect of the strength of measurement noise is studied in Sec.~\ref{sec:noise} of the SM.

In cases where the system itself is not fully deterministic, Eq.~\eqref{eq:x1} becomes a set of SDEs. Opsin models are population models (i.e., with conservative dynamics), and their ODEs can be extended to the following SDEs

\begin{align}
{\rm d }{\bm X} &= \bm F(\bm X,\, u(t);\, {\bm \Theta})\,{\rm d}t + {\bm G}(\bm X)\,{\rm d}\bm W\,,\label{eq:sde}
\end{align}
where $\bm W$ is a vector of Brownian processes and $\bm G$ is the function combining them into the noise signal. 
We chose $\bm G(\bm X) = \alpha(\mathbb{I}_n-\bm X\bm 1^\top)\sqrt{\bm X}$ as it enforces that $0\le x_i\le 1$ and $\sum_ix_i=1$ as required. Fig.~\ref{fig:noise} (II) shows that even with a noisy reference system, the scheme still selects the correct model. The stochastic nature of the system tends to prevent the parameter fitting process from overfitting.
Note that the stochastic terms are added only to the reference system; we still treat the candidate models as being noise-free during both the fitting and discrimination processes.

Our numerical investigation shows our scheme is capable of discriminating models even when the models are only idealizations of the true system. The scheme being validated, we can now move to a lab setting in which the ground truth process is not known.

\section{Lab Experiments}\label{sec:exp}
Finally, we test our algorithm in a closed-loop using an experimental system involving mammalian cells transfected to express the ChR2 channelrhodopsin protein.
Recordings are obtained via voltage clamp in a whole cell patch clamp configuration as illustrated in Fig.~\ref{fig:patch}.
For brevity, we here outline only the key details, with full methods provided in Sec.~\ref{sec:expt_methods} of the SM.
The voltage clamp techniques holds the potential difference across the cell membrane fixed at a specified value (here -60 mV) and measures currents flowing across the cell membrane.
By using the control signal $u(t)$ as an input to a 470 nm LED, as illustrated in Fig.~\ref{fig:patch}, we can thus measure the photocurrents induced as the opsin proteins change state.

\subsection{Hardware-Software Setup}
To avoid issues associated with physiological run down (i.e., slow loss of electrophysiological signal) and loss of overall cell viability, it is imperative to conduct patch clamp recordings in as short a time as is reasonably possible.
The duration of open loop recordings is typically known and set in advance to ameloriate such problems.
By contrast, our scheme requires solving both a minimization and an optimal control problem at each iteration.
Thus, whilst the maximum possible \emph{recording} time, that is, the time in which photocurrents are being measured, is known in advance, the total real-world time required for the algorithm to complete is not.
It is thus important the parameter fitting and control input optimization iterations to complete as quickly as possible, and should ideally take less time than the recordings themselves.

To achieve high computational efficiency, particularly, accounting for older computational hardware in the lab, we leverage the strength of Symphony, a Matlab-based data acquisition package, Transmission Control Protocol (TCP) servers, and the Julia programming language to accelerate the computations required for our system, as illustrated in the right-most portion of  Fig.~\ref{fig:patch}.
In this configuration, the computational part of the method is implemented in a Julia code that is precompiled and called only when data have been collected.
Conversely, Symphony waits until the next control signal has been computed by the Julia server and loads the resulting profile to send to the LED.
Further computational speed gains can be achieved by noting that the Julia server need not be running on the same computer as the one collecting the data.
We thus take advantage of high speed intranet to offload the computational work to a dedicated computational server.
The latter had the following CPUs: 12th Gen Intel(R) Core(TM) i7-12700\,.

\subsection{Experimental Results}
Instead of using the termination criterion Alg\,.~\ref{alg:stop}\,, we rather run the experiment for a large fixed number of iterations. 
We do this to collect as much data as possible, keeping in mind that the patched cell might be lost during the experiment.
Nevertheless, as we will see below, even when the cell is lost during the experiment, one may obtain useful information from our model discrimination algorithm.

\begin{figure*}[t!]
    \centering
    \includegraphics[width=\textwidth]{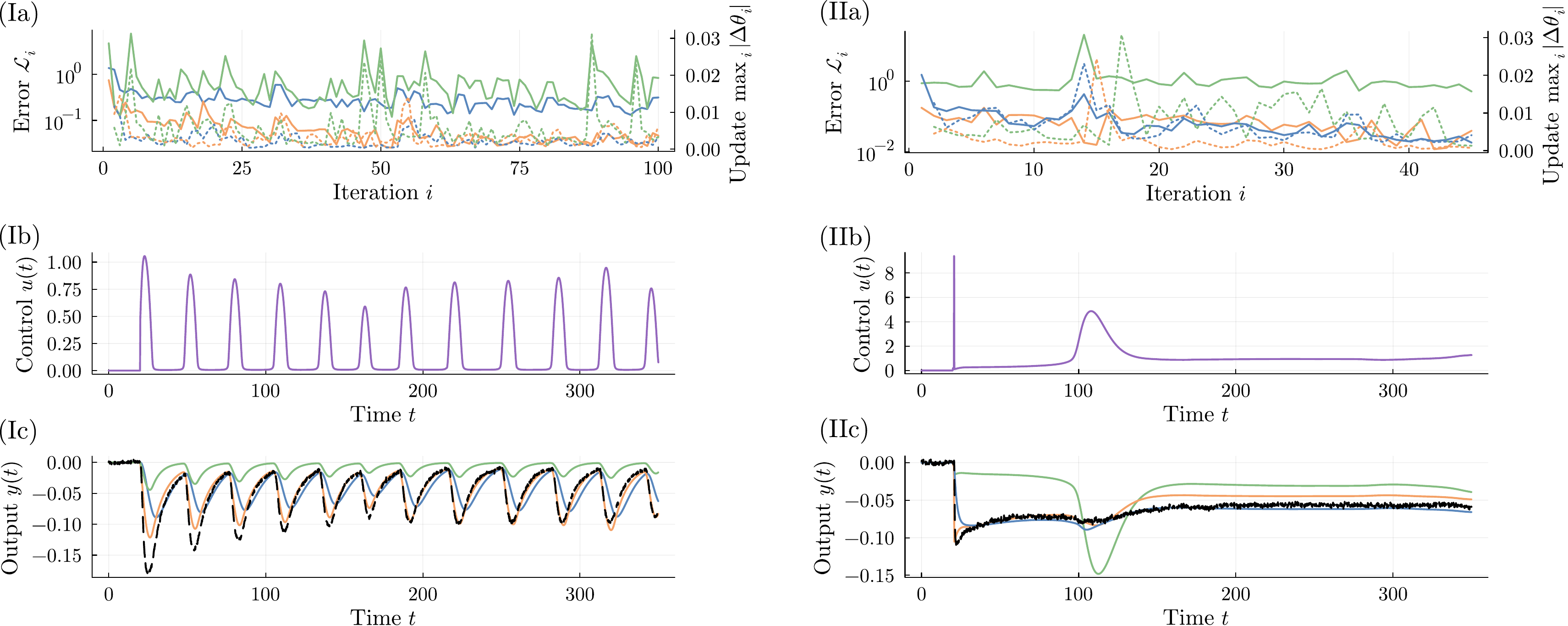}
    \caption{Two laboratory experiments, labeled I and II, discriminating 3-state (green), 4-state (orange) and 6-state (blue) models. (a) Evolutions of mismatch with the reference trajectory (solid) and of maximal parameter increment (dashed). (b) Optimal control at the end of the procedure. (c) Corresponding responses of the models and of the experiment (black).}
    \label{fig:exp}
\end{figure*}

In the first experiment, we compare the 6-state and 4-state models and then the 3-state and 4-state models. For the first comparison, we were able to run 100 iteration of the algorithm, but for the second, we lost the cell after iteration 46.
In Fig.~\ref{fig:exp} (I), we discriminate between the 6-state and 4-state model and evaluate the 3-state model. The top panel shows that the error decays faster for the 4-state (red line) than for the 6-state (blue line) and 3-state (green line) models. For both models, the parameters are still fluctuating between iterations, with seemingly smaller variation for the 4-state model (see dashed lines in top panel).
In Fig.~\ref{fig:exp} (Ib) and (Ic), we respectively show the final control input designed and the response of the experiment and the models. 
We emphasize here that the blue, orange and green lines are not obtained by performing a fit to the data, but from predicting the response using the mathematical models and the parameters inferred at the previous iteration. 
This designed control input seems efficient in highlighting the different response mechanisms of the two models. 
Indeed, the 4-state model include a direct transition from the closed to the open states enabling a swift reaction to fast varying light inputs unlike the 6-state model, for which a transition to intermediary state is necessary.
The 4-state model remarkably predicts the response of the experiment.
Therefore, based on this model discrimination result, one would choose the 4-state over the 6-state model to predict the future evolution of the experiment. 
For reference, the parameter estimates at each iteration are shown in the SM in Fig.~\ref{fig:case_11} with representative control inputs plotted in Fig.~\ref{fig:case_11b}.

Next, we compare the 3-state and 4-state models in Fig.~\ref{fig:exp} (II). 
The algorithm ran until iteration 47, after which, we lost the cell.
In Fig.~\ref{fig:exp} (IIa)l, one observes that the error in the model prediction for the 4-state model (orange line) decays while for the 3-state (green line), it remains approximately constant. One also observes that the parameters for the 3-state model are still varying considerably compared to the 4-state model. With these control signals, the 6-state model (blue line) performs as well as the 4-state model. 
In Fig.~\ref{fig:exp} (IIb) and (IIc)\,, we respectively show the last control input that was designed and the actual response in the experiment (gray line) and the predicted response from the 3-state (green) and 4-state (orange) models.
The 4-state outperforms the 3-state in predicting the outcome of the experiment.

\section{Discussion and Future Work}\label{sec:discussion}
The problem of predicting the future evolution of an observed system is of primal importance in many fields. This task is achievable by selecting and analyzing a mathematical model. 
Here, we propose a closed-loop framework to pick the most adequate mathematical model over a set of candidates when the observed system can be controlled. 
At each iteration of our method, we sequentially perform parameter fitting of the candidate models to the data, and then devise control input signals that maximally differentiate the output of the models. 
The generated control signal is eventually re-injected into the observed system that will be used for the next iteration. 
We here focused specifically on opsin models, due to the ongoing interest they attract in experiments using optogenetics. 
First, our results highlight that using control inputs that are not complex enough do not allow to perform model discrimination. 
Indeed, using a box-shaped control signal, all the competing models well reproduced the data. 
Second, optimizing a performance metric that both quantifies the difference between the trajectories of the candidate models and help parameter fitting over the iterations allows to select the model that better predicts the observed system. 

Finally, we apply the method precisely in the setting for which it has been designed, namely, for choosing over a set of candidate models, the one that is the most accurate to predict the response of an actual experiment. 
We performed an electrophysiology experiment where the control input is a light signal, and the output is the resulting photocurrent induced by the stimulation of channelrhodopsin. 
To predict the photocurrent, we considered three candidate models for the opsin dynamics that include different mechanisms. 
We built an efficient hardware-software interface that enables the online application of our closed-loop algorithm. 
Using this setup, our method identified the 4-state model as the most accurate to predict the photocurrent in our experiment. 
Importantly, this demonstrates the practicability of our method to discriminate candidate models on real-world settings.

When discriminating between more than two prospective, our framework can be applied in a binary tournament format in which candidates are compared pairwise, with the winning model compared to other winning models until a sole winner is found.
In principle, it may be possible to derive a performance metric ($\mathcal{J})$, for example using the variance of the model responses over time (see Sec.~\ref{SimModDisc} in the SM), that facilitates optimal discrimination between more than two models simultaneously.
Alternatively, one could use a Fokker--Planck formulation of \eqref{eq:sde} and then minimize the overlap between the output distributions for the observed variable, noting that different models likely filter noise in different ways.
However, the requirement to solve PDEs for each model, rather than ODEs would lead to importantly increased computational costs.
Overall, we leave the development of routines to discriminate between more than two models simultaneously as an avenue for future research.

While we demonstrated the effectiveness of our method in numerical simulations and on an electrophysiology experiment, we do not provide a mathematical proof that the algorithm will generically converge to the most accurate candidate model. Indeed, such proof is challenging to obtain if possible at all for real-world data. In Sec.~\ref{app:conv} of the SM, we provide evidence of local convergence of the algorithm under restricted conditions. 

Our framework has been applied to an electrophysiology experiment involving optogenetic stimulation.
The use of light-based stimulation allowed us to isolate the dynamics of one ion channel type, which greatly facilitates model discrimination. 
In a more general setting, one could adapt our model discrimination framework to current or voltage clamp data from excitable cells under electrical stimulation.
This setting requires discriminating model descriptions of multiple channel types simultaneously, taking advantage of the fact that different ion channel gates typically operate over different timescales.
Previous studies have addressed this problem by generating a large cell model containing models for all ion channel types that may be present in the cell and then by designing current inputs to optimize parameter fitting~\cite{nogaret_automatic_2016} using long time-series recordings.
During the fitting process, the conductance parameter of ion channels (which scales the contribution of that part of the model to the measured output) not present in the cell converge to zero, and hence these sub-components of the model can be removed.
This approach produced models that successfully predicted the impact of pharmacological agents on the cell's response, validating the model selection procedure. 
However, the process of fitting the larger model imposed an important computational cost on the fitting procedure, leading to fitting processes taking several hours, meaning that fitted models could not be tested using the same cell for which the model was built.
In addition, this same cost prohibited the comparison of different model structures for each of the ion channel types as we do in this manuscript.
We anticipate that using an iterative scheme involving shorter recordings may help address both of these issues, leading to improve model discovery, particularly, if the process can be combined with genetic information about which proteins may be expressed in the cell~\cite{Lipovsek2021}.

The above requires applying our approach to a potentially very high-dimensional model, compared with the relatively low-dimensional systems we here consider.
In our experiments, the major computational bottleneck was in the optimization of the control signal, primarily as a result of our choice to use IPOPT to solve \eqref{opt}-\eqref{eq:bound}.
For larger models, the optimal control input can be obtained efficiently using numerical techniques such as the adjoint method, as used in \cite{liu2025optimal} whose sensitivity to system size is limited.
In this case, the major computational bottleneck arises in the parameter fitting step, which is a challenge faced in most computational biology models.
Here, we have some control over the complexity through the choice of hyperparameters, such as the learning rate of the Adam algorithm and the number of gradient descent steps taken in each iteration of our loop.
Balancing the per iteration time with the number of iterations required for convergence likely requires empirical experimentation on a system-by-system basis.
We remark, however, that the scheme is intended to be used in an online fashion during a wet lab experiment and so models with large numbers of state variables and parameters may be inappropriate unless many of the parameters are specified in advance.

Real-world systems, particularly in biology, exhibit large heterogeneity and are influenced heavily by noise.
From a modeling perspective, this leads to multiple sources of uncertainty.
These sources could be incorporated into models in a variety of ways, including through additional terms to tackle model discrepancy, and by making estimates for parameters uncertainty.
The latter of these could, for example, be achieved by switching to a Bayesian framework \cite{chalonerBayesianDesign}, in which parameters are described by distributions rather than by point estimates. 
Although this switch considerably increases the computational workload, since we would need to estimate posterior distributions, the Bayesian framework naturally captures the variability in responses from the cell and provides a means to naturally incorporate measurement noise.
Moreover, through the use of hierarchical Bayesian modeling, we could assess variability across populations of cells, which may be useful when when considering the impact of diseases or genetic mutations, since these may lead to shifts in distributions of parameters that cannot be assessed through point estimates alone~\cite{saghafi_inferring_2024}.

One of the main limitations with our current approach is that both the parameter fitting and control design use gradient information in their optimization.
As such, the performance of our scheme is sensitive to our initial estimates for parameters and thus we cannot guarantee that we find global, rather than local, optimal parameter values and control inputs even though the parameters of each of the Markov models are structurally identifiable (see Sec~\ref{app:id} of the SM).
For similar reasons, we also cannot guarantee that we do not overfit the data, despite adding penalty terms to the control optimizer to force the system to explore different parts of state space.
As such, a post processing may be required to search for the globally optimal parameter values.
One potential solution to this problem that could be used during the algorithm run is to use a surrogate model, such as a Gaussian Process, to approximate the model in locations not sampled~\cite{mohammadi_cross-validation--based_2022}.
In this sense, a more global search for optimal parameters and control inputs can be performed in which the control inputs balance between model discrimination and improving the surrogate model through adaptive sampling.
The development of surrogate models could also be useful in amortization schemes, which front load computations prior to any experiment \cite{fosterDeepAdaptiveDesign,huang2025aline}, since they would allow for policies for control design across a distribution of parameter combinations, circumventing the requirement to have strong prior parameter assumptions.
The drawback of this approach is that building a ``good enough'' surrogate model likely requires a large amount of initial data collection, so careful consideration of the overall experiment duration needs to be made.

The structure of the three Markov models we contrast is assumed to be known, so that only parameters of the models need to be fit in the first part of each algorithm iteration.
A more general approach would involve discovering the equations using techniques from Scientific Machine Learning and Symbolic Regression, such as SINDy.
In this case, a targeted library of possible functions or a basis set of functions and compositions could be provided to build a model `from the ground up'.
Importantly, such basis functions would be selected based on their relevance and plausibility for the system of interest, rather than using standard libraries, such as polynomials up to certain degree.
In this way, the combination of our approach with model discovery algorithms could be used to build interpretable models in an online fashion.
Alternatively, since symbolic regression approaches typically output multiple candidate models (either as ranked list for a single set of rules, or through specification of different allowable interactions between terms in different regression runs), our approach could be used as additional step to determine which of these is most appropriate.
Control targets could then be designed to maximize the information gained regarding whether to incorporate specific terms, as discussed in~\cite{bao_information_2025}.
Such an approach could provide a flexible means to build interpretable models from partially observed systems that avoid existing bias in the literature.
However, since the identified functional forms may focus on particular parts of the time series, extra care must be taken to avoid overfitting.

Another potential limitation of our scheme is in parameter inference in certain cases. Even in cases where parameters are structural identifiable, if some of them are close to be unidentifiable, for example if the output only depends on them in a distant and weak way, we expected our method to struggle inferring their values. Similarly, if the model’s sensitivity to a parameter is only non-zero for a peculiar control, our scheme might never gain information about that parameter. This could be addressed for instance by adding a random control profile every now and then. Going beyond identifiability, if the goal is to obtain the best predictive model so that system can be controlled efficiently, the model selected by our method will perform as long as expected control is similar to control used during the selection process, which can be tested.

Since our model discrimination framework requires only the possibility to observe time series output from and apply time series control to a given system, we anticipate that it can easily be adapted beyond biological systems.
For example, one could investigate which model is the most adequate to predict other dynamics, such as opinion formation on influence networks, trades on financial markets, essentially using the algorithm for developing digital twins.


\section*{Acknowledgments}
The research was funded by the UKRI grant number MR/X034240/1\,.

\section*{Declaration of interest}
The authors declare no competing interests.

\section*{Author Contributions}
Conceptualization: LP, MT, KW. 
Formal Analysis: LP, MT.
Resources: AJ, PT, KW.
Visualization: LP, MT, KW.
Methodology: LP, MT, KW.
Writing - original draft: LP, MT, KW.
Writing - review and editing: LP, MT, KW.
Funding acquisition: KW.

\section{Data and Materials Availability}
All data needed to evaluate the conclusions in the paper are present in the paper, the Supplementary Materials and the GitHub repository~\cite{github} where we uploaded the code implementing the presented algorithms.


\pagebreak
\pagebreak

\widetext
\begin{center}
\textbf{\large Supplementary Materials: A Closed-loop Framework to Discriminate Models Using Optimal Control}
\end{center}
\setcounter{equation}{0}
\setcounter{figure}{0}
\setcounter{table}{0}
\setcounter{page}{1}
\makeatletter
\renewcommand{\theequation}{S\arabic{equation}}
\renewcommand{\thefigure}{S\arabic{figure}}
\section{Description of Advanced Opsin Models}\label{opsins}

A slightly more advanced description of opsin dynamics compared to 3-state model \eqref{3s1}-\eqref{3s2} is given by the following  4-state model, which was introduced in~\cite{hegemannMultiplePhotocyclesChannelrhodopsin2005} to account for the biexponential decay of photocurrents and the observation of light-dependent deactivation and recovery rates
\begin{subequations}
\begin{align}
\dot x_1 &= \theta_1  u(t)\, x_4 + \theta_2 u(t)\, x_2 - \big[\theta_3 +\theta_4 u(t) \big] x_1\,,\label{4s1}\\
\dot x_2 &= \theta_5  u(t)\,  x_3 + \theta_4 u(t)\, x_1 - \big[\theta_6 +\theta_2 u(t) \big] x_2\,,\\
\dot x_3 &= \theta_6 \, x_2 - \big[\theta_7 + \theta_5 u(t) \big] x_3\,,\\
x_4 &= 1 - x_1 - x_2 - x_3\,,\\
y & = \theta_9\, x_1 + \theta_8\, x_2\,. \label{4s2} 
\end{align}
\end{subequations}
This 4-state model possesses two open ($x_1$ and $x_2)$ and two closed states ($x_3$ and $x_4$) with the observed photocurrent \eqref{4s2} taken to be proportional to the weighted sum of the proportions of opsins in the open state.

A yet more advanced model is presented in~\cite{grossmanSpatialPatternLight2013}, which introduces two additional states, referred to as \textit{intermediate states}, which account for the second order dynamics in the transitions from the closed to the open states.
The governing equations for this 6-state model are
\begin{subequations}
\begin{align}
\dot x_1 &= \theta_1  u(t)\, x_6 - \theta_2  x_1\,,\\
\dot x_2 &= \theta_2 x_1 + \theta_3 u(t)\, x_3 - \big[\theta_4 +\theta_5 u(t) \big] x_2\,,\\
\dot x_3 &= \theta_6  x_4 + \theta_5 u(t)\, x_2 - \big[\theta_7 +\theta_3 u(t) ] x_3\,,\\
\dot x_4 &= \theta_8 u(t)\,  x_5- \theta_6 x_4\,,\\
\dot x_5 &= \theta_7 \, x_3- \big[\theta_9 +\theta_8 u(t) \big] x_5\,,\\
x_6 &= 1 - x_1 - x_2 - x_3 - x_4 - x_5\,,\\
y & = \theta_{11}\,x_2 + \theta_{10}\, x_3\,,
\end{align}
\end{subequations}
where $x_2$ and $x_3$ are open states, $x_5$ and $x_6$ are closed states, and $x_1$ and $x_4$ are intermediate states.

To account for the desensitization effect often observed in channelrhodopsins under prolonged or repeated light stimulation, \cite{FrommeThesis} postulated an 8-state Markov model, the spectroscopic signatures of which were later observed experimentally~\cite{kuhne_unifying_2019}.
\begin{subequations}
\begin{align}
\dot x_1 &= \theta_1 u(t) \, x_8 - \theta_2 \, x_1, \\
\dot x_2 &= \theta_2 \, x_1 - \left[ \theta_3 + \theta_4 u(t) \right] x_2, \\
\dot x_3 &= \theta_4 u(t) \, x_2 + \theta_5 \, x_6 - \theta_6 \, x_3, \\
\dot x_4 &= \theta_6 \, x_3 - \left[ \theta_7 + \theta_8 \right] x_4, \\
\dot x_5 &= \theta_7 \, x_4 + \theta_3 \, x_2 - \theta_9 \, x_5, \\
\dot x_6 &= \theta_{11} u(t) \, x_7 - \theta_5 \, x_6, \\
\dot x_7 &= \theta_8 \, x_4 - \left[ \theta_{10} + \theta_{11} u(t)\right] x_7, \\
x_8 &= 1 - x_1 - x_2 - x_3 - x_4 - x_5 -x_6 - x_7 \,,\\
y & = \theta_{12}\,x_2 + \theta_{13}\, x_3\,.
\end{align}
\end{subequations}

\begin{figure}[h!]
\includegraphics[width=0.4\columnwidth]{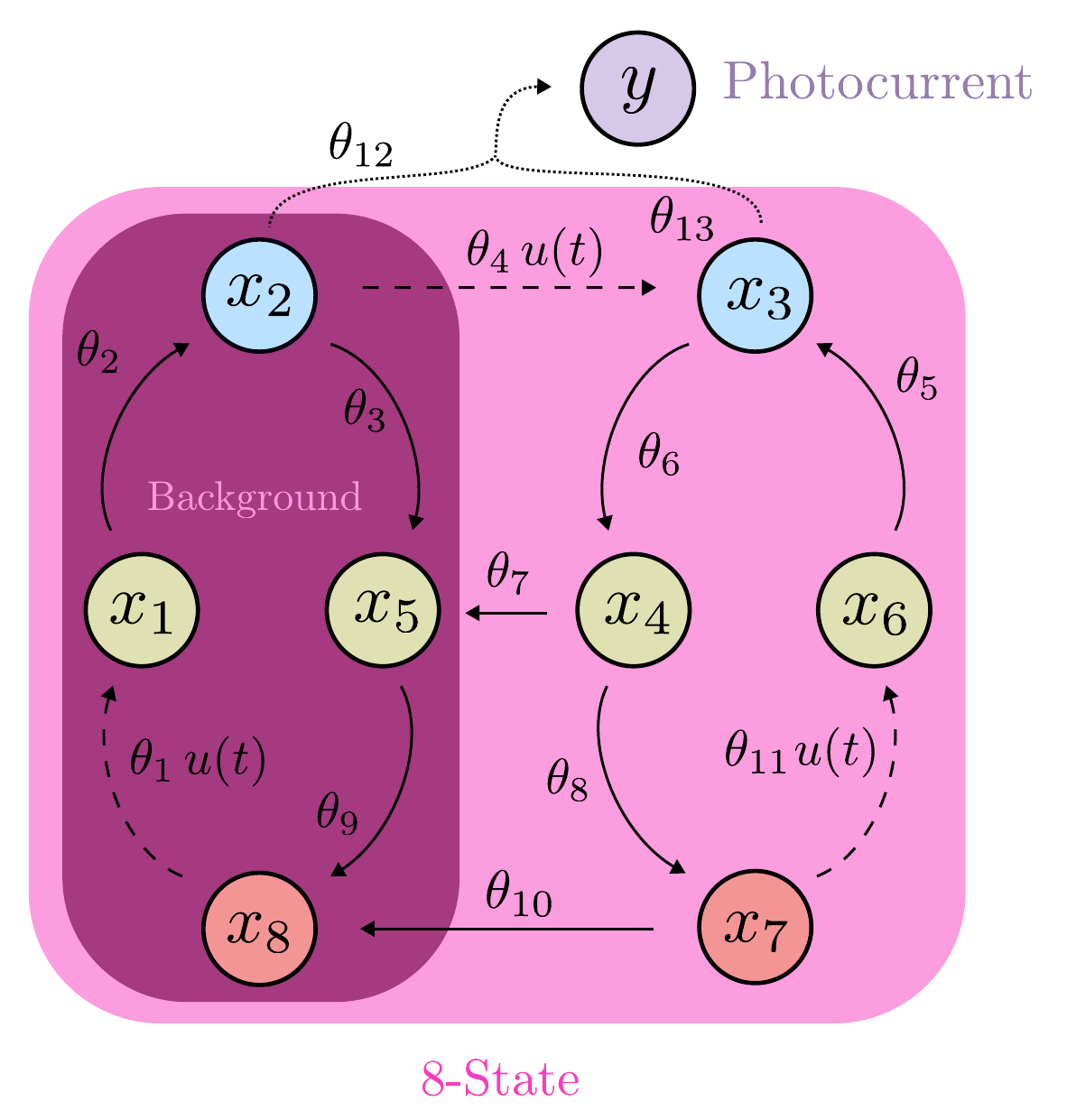}
\caption{Depiction of 8-state model. Color convention is the same as in Fig.~\ref{fig:opsin}.}\label{fig:8state}
\end{figure}

The model splits the opsin light cycle into two sub-cycles, one corresponding to a background (with states $x_1$, $x_2$, $x_5$, $x_8$) in which the opsin has not been exposed to light for an extended period of time, the other corresponding to adaptation under light exposure (with states $x_3$, $x_4$, $x_6$, $x_7)$.
Each sub-cycle comprises a ground state, an intermediate state, an open state, and a closed state.
Under exposure to light, an opsin in a ground state ($x_7$ or $x_8$) moves to an open state ($x_2$ or $x_3$) through an intermediate state ($x_1$ or $x_6$).
An opsin in the background open state, $x_2$, in may then either move to the background closed state, $x_5$, or to the adapted open state, $x_3$ depending on the light intensity.
Similarly, an opsin in the adapted closed state  $x_4$ may transition to the adapted ground state $x_7$ or to the background closed state, $x_5$.
Finally, closed states ($x_4$, $x_5$) transition to ground states, and the adapted ground state may also transition to the background ground state, $x_8$. The 8-state model is depicted in Fig.~\ref{fig:8state}.

\begin{table}[h!]
\begin{tabular}{lrclrlrclrlr}
\hline
\multicolumn{2}{c}{3-State}	&\hphantom{aa}&	\multicolumn{4}{c}{4-State} 					&\hphantom{aa}&			\multicolumn{4}{c}{6-State} 							\\		
\hline																			
$\theta_1$:	&	0.24	&&	$\theta_1$:	&	0.250,	&	$\theta_6$:	&	0.089	&&	$\theta_1$:	&	1.404,	&	$\theta_7$:	&	0.065	\\
$\theta_2$:	&	0.053	&&	$\theta_2$:	&	0.134,	&	$\theta_7$:	&	0.014	&&	$\theta_2$:	&	0.883,	&	$\theta_8$:	&	0.003	\\
$\theta_3$:	&	0.02	&&	$\theta_3$:	&	0.130,	&	$\theta_8$:	&	-0.235	&&	$\theta_3$:	&	0.077,	&	$\theta_9$:	&	0.032	\\
$\theta_4$:	&	-0.533	&&	$\theta_4$:	&	0.116,	&	$\theta_9$:	&	-0.27	&&	$\theta_4$:	&	-0.001,	&	$\theta_{10}$:	&	-0.825	\\
	&		&&	$\theta_5$:	&	0.004	&		&		&&	$\theta_5$:	&	0.644,	&	$\theta_{11}$:	&	-0.186	\\
	&		&&		&		&		&		&&	$\theta_6$:	&	0.804	&		&		\\
\hline
\end{tabular}

    \caption{3-, 4- and 6-state models fitted to experimental data used as ground truth (reference) models.}
    \label{tab:ref_param}
\end{table}

\section{Model discrimination with the 8-state model}
\begin{figure}
    \centering
    \includegraphics[width=0.5\linewidth]{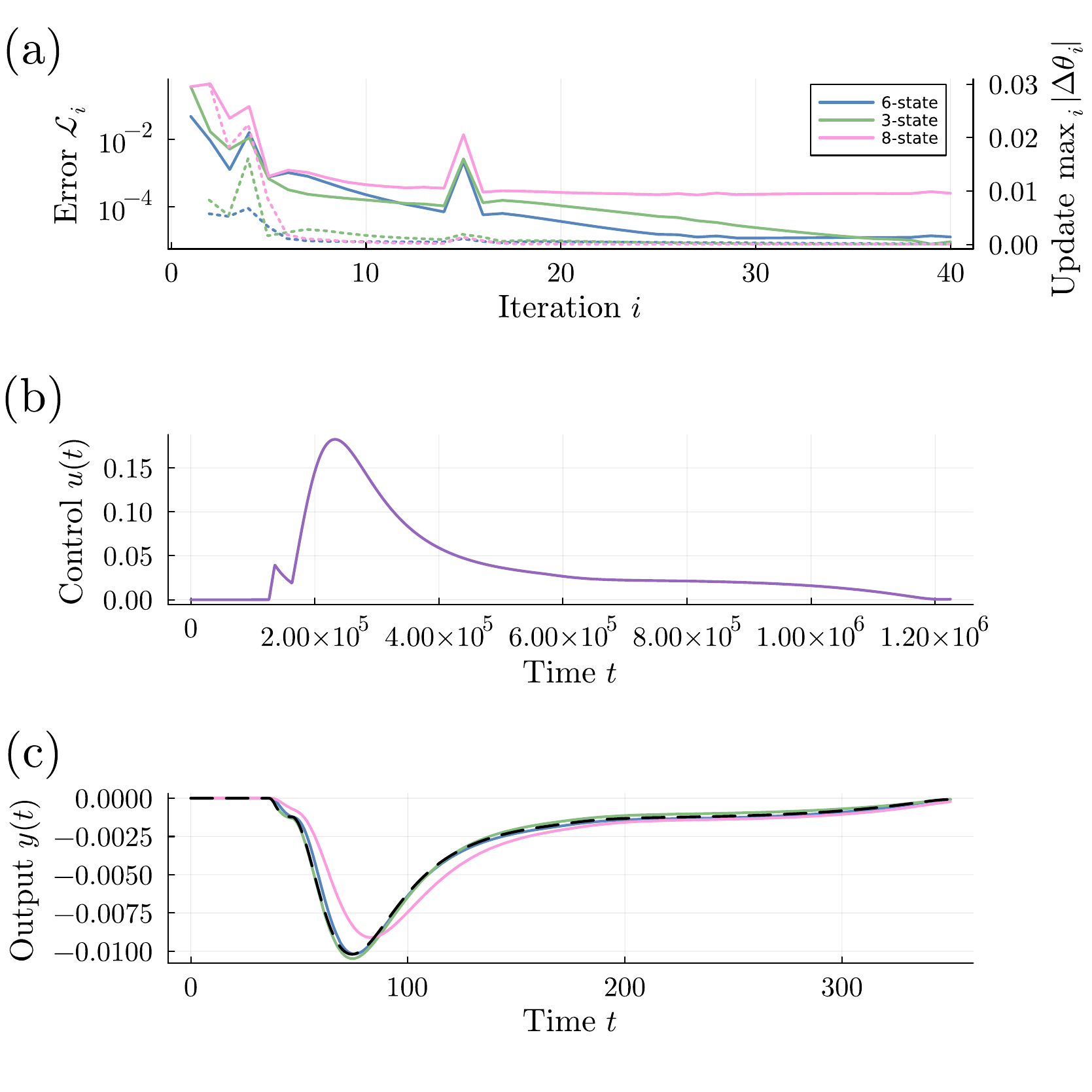}
    \caption{(a) Evolution of mismatch with the reference trajectory (solid) and of maximal parameter increment (dotted) for
the 3-state (green), 6-state (blue) and 8-state (pink) models. The control input are obtained by discriminating between the 3- and 8-state models. The parameters of the 6-state model are fitted and the error in the prediction is plotted, but it does not enter the model discrimination loop. (b) Optimal
control at the end of the procedure. (c) Corresponding responses of both models compared with the reference model
(dashed black), which in this case is the 3-state model with
tuned parameters.}
    \label{fig:8-state_ref3}
\end{figure}

\begin{figure}
    \centering
    \includegraphics[width=0.5\linewidth]{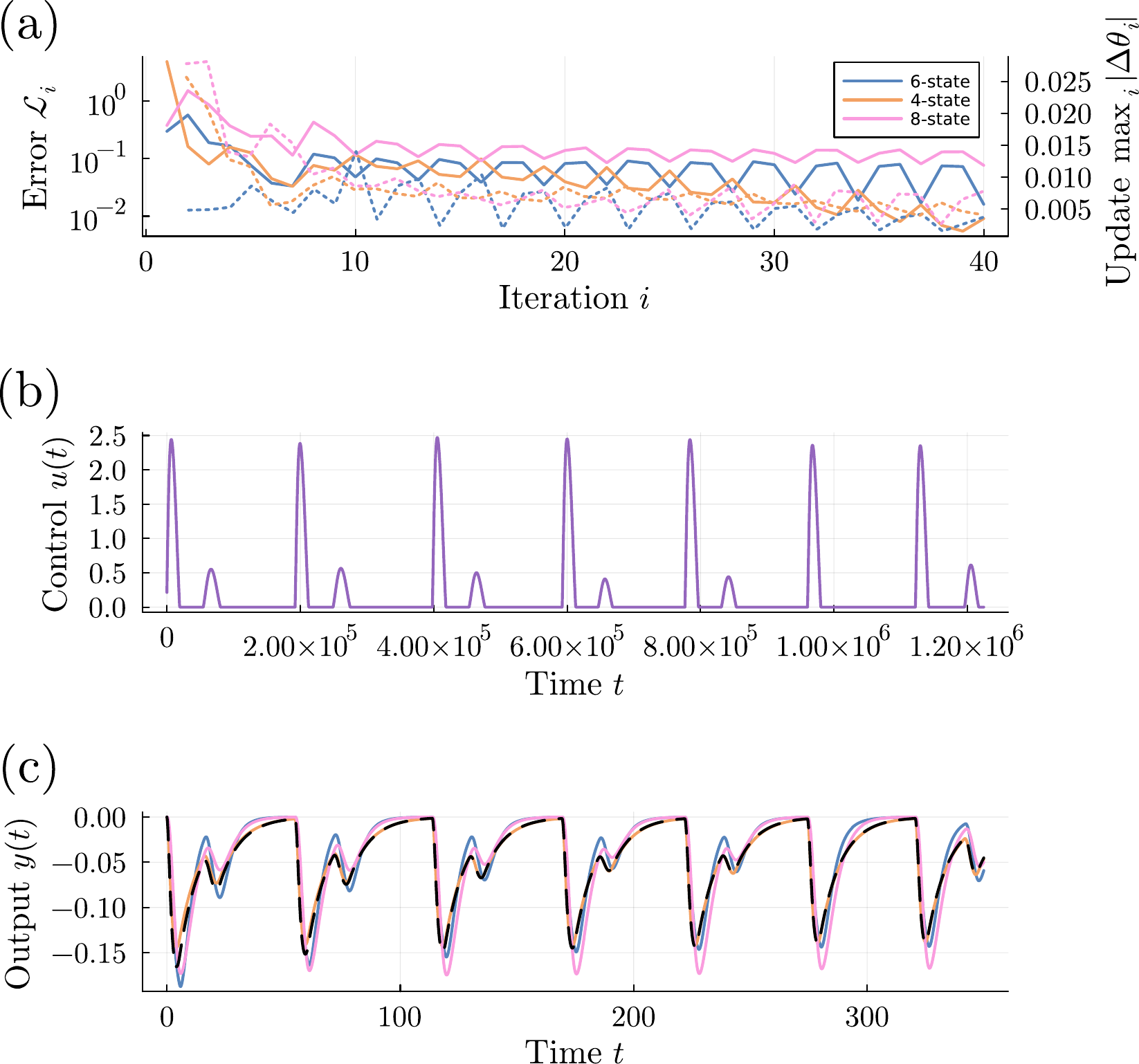}
    \caption{(a) Evolution of mismatch with the reference trajectory (solid) and of maximal parameter increment (dotted) for
the 4-state (orange), 6-state (blue) and 8-state (pink) models. The control input are obtained by discriminating between the 4- and 8-state models. The parameters of the 6-state model are fitted and the error in the prediction is plotted, but it does not enter the model discrimination loop. (b) Optimal
control at the end of the procedure. (c) Corresponding responses of both models compared with the reference model
(dashed black), which in this case is the 4-state model with
tuned parameters.}
    \label{fig:8-state_ref4}
\end{figure}

We perform model discrimination between the 8-state and the 3-state models while the reference is obtained from a 3-state model, and the 8-state and the 4-state models, while the reference is obtained from a 4-state model. The results are shown in Figs.~\ref{fig:8-state_ref3} and \ref{fig:8-state_ref4}, where one obtains the correct identification of the 3-state and 4-state models respectively. 

\section{Performance measure beyond pairwise model discrimination}\label{SimModDisc}
\begin{figure}
    \centering
    \includegraphics[width=0.5\linewidth]{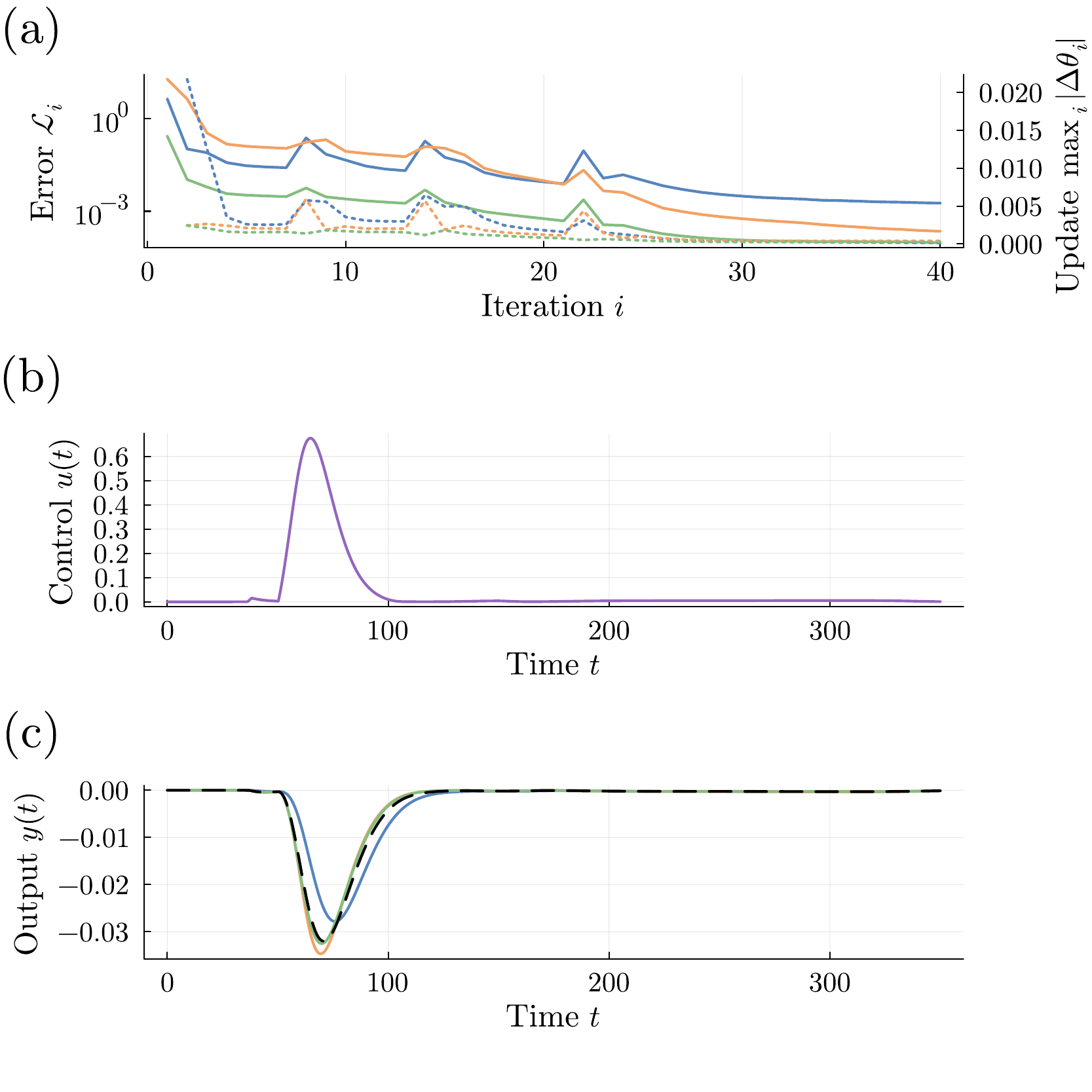}
    \caption{a) Evolution of mismatch with the reference trajectory (solid) and of maximal parameter increment (dotted) for
the 3-state (green), 4-state (orange) and 6-state (blue) models. The control input are obtained by discriminating between the three models using Eq.~(\ref{cost_var}). (b) Optimal
control at the end of the procedure. (c) Corresponding responses of both models compared with the reference model
(dashed black), which in this case is the 3-state model with
tuned parameters.}
    \label{fig:cost_var3}
\end{figure}

\begin{figure}
    \centering
    \includegraphics[width=0.5\linewidth]{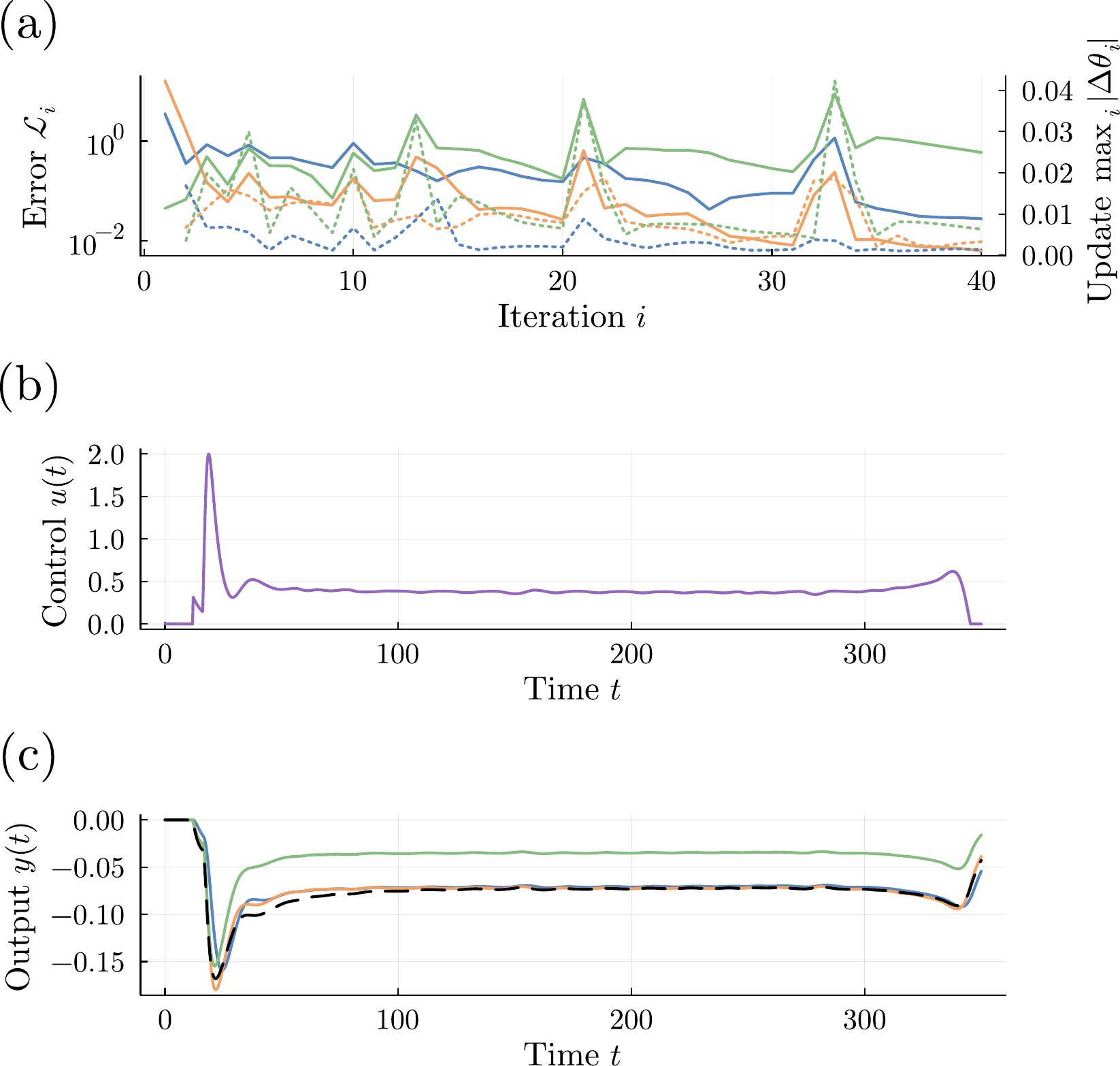}
    \caption{(a) Evolution of mismatch with the reference trajectory (solid) and of maximal parameter increment (dotted) for
the 3-state (green), 4-state (orange) and 6-state (blue) models. The control input are obtained by discriminating between the three models using Eq.~(\ref{cost_var}). (b) Optimal
control at the end of the procedure. (c) Corresponding responses of both models compared with the reference model
(dashed black), which in this case is the 4-state model with
tuned parameters.}
    \label{fig:cost_var4}
\end{figure}

The algorithm discussed in the main text performs model discrimination by comparing pairs of candidate models. 
Instead of comparing pairs, one could compare more than two models at a time. 
A possible performance metric to achieve such a discrimination, corresponds to maximizing the variance of the candidate model responses. Here, we assume that the output are 1-dimensional time-series, which is the case of the opsin models we considered.
The metric reads,
\begin{align}\label{cost_var}
u^* &= {\rm argmax}_{u}\,\mathcal{J}_{\rm mult}(u)\,,\\
    \mathcal{J}_{\rm mult}\big(u(t)\big) &:= \int_0^T \Big[{\rm Var}\Big(\{ y_i\big(\bm X_i(t)\big)\}_{i=1}^{N_c}\Big)- C(u)\Big]\,{\rm d}t\,, \nonumber
\end{align}
with $N_c$ the number of candidate models and 
\begin{align}
    {\rm Var}\Big(\{ y_i\big(\bm X_i(t)\big)\}_{i=1}^{N_c}\Big) = N_c^{-1}\sum_{i=1}^{N_c}\left(y_i(t) - \overline{y}(t)  \right)^2\,,
\end{align}
where $\overline{y}(t) = N_c^{-1}\sum_{i=1}^{N_c}y_i(t)$ and $C$ defined as in the main text. 
In Figs.~\ref{fig:cost_var3}\,, \ref{fig:cost_var4}\,, we perform the model discrimination among the 3-, 4- and 6-state models simultaneously using the loss function (\ref{cost_var})\,, when the reference model is respectively, a 3-state and a 4-state model. 
While this approach allows to discriminate all the models simultaneously, the optimization of the Eq.~(\ref{cost_var}) takes at least an order of magnitude more time than optimizing the pairwise loss function of the main text. 
That makes it impractical for the experiment we are considering in this work. However, such multi-model discrimination might be valuable for experiments in which the time between iteration is less critical.

\section{Effect of Measurement Noise}\label{sec:noise}

\begin{figure}
    \centering
    \includegraphics[width=0.5\linewidth]{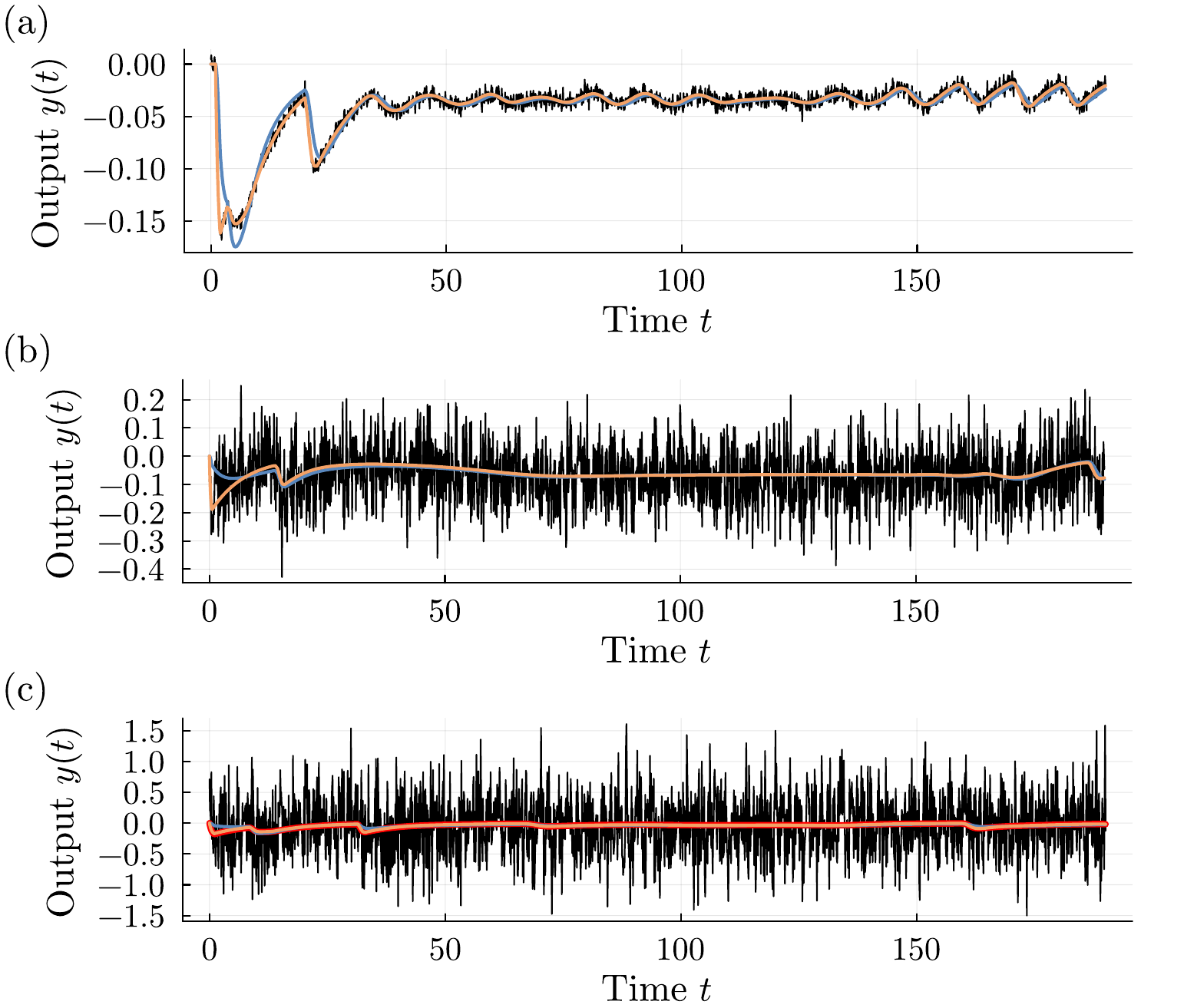}
    \caption{Response of the 4-state (orange) and 6-state (blue) models against the reference model trajectory (black), here the 4-state model, with different levels of noise $\xi=0.005,0.1 \text{ and }0.5$ in panels (a), (b) and (c) respectively. On panel (c), the red line displays the reference trajectory without noise.}
    \label{fig:noise_to_signal}
\end{figure}

Fig.~\ref{fig:noise_to_signal} studies when noise breaks our scheme. It shows that our scheme is surprisingly robust to noise. Our take on this matter is that for white noise, the uncertain signal is unbiased: it quickly oscillates between overestimating the noiseless signal and underestimating it. The learning step cannot reach the target, and the error remains high, but it still minimizes the mismatch with the noisy reference, thus learn the mean trajectory which is equivalent to the reference trajectory. In case of time-correlated noise, e.g., Ornstein–Uhlenbeck process, the situation would be different and we expect the scheme to fail  for lower levels of noise.

\section{Experimental methods}
\label{sec:expt_methods}

\paragraph*{Cell culture}
The GH3 epithelial-like rat pituitary tumor cell line (ATCC CCL-82.1, American Type Culture Collection, USA) was maintained in Ham's F-12 Nutrient Mixture (Sigma-Aldrich, USA) at \SI{37}{\degree C} in a humidified atmosphere containing 5\% CO$_2$. Culture medium was supplemented with 15\% horse serum, 2.5\% fetal bovine serum, and 1\% GlutaMAX (Gibco). Cells were routinely passaged at $\sim$80\% confluency using a 1:3 subculture ratio. For experiments, cells were seeded into 12-well plates at a density of 1 $\times$ 106 cells per well and allowed to stabilize for 24h before transfection.

\paragraph*{Transfection} The channelrhodopsin-2 (ChR2) expression plasmid pRP[Exp]-EGFP-CAG$>$ChR2 was obtained from VectorBuilder (USA). GH3 cells were transfected using TransIT-2020 Transfection Reagent (Mirus Bio, USA) at a 3:1 reagent:DNA ratio. TransIT-2020:DNA complexes were prepared in serum-free growth medium (SFM) by combining 1.2 µg of plasmid DNA (200 ng/µL), 3.6 µL TransIT-2020 reagent, and 100 µL SFM per well, followed by incubation at room temperature for 15 min. The resulting complexes were added dropwise to each well, and cells were incubated at \SI{37}{\degree C} / 5\% CO$_2$ for 48h prior to medium change for subsequent electrophysiology experiments.

\paragraph*{Patching solutions}
The extracellular solution comprised 2.8 mM glucose, 132 mM NaCl, 5 mM KCl, 2.6
mM CaCl2, 1.2 mM MgCl2, 10 mM HEPES (pH 7.4 adjusted with NaOH,
osmolarity 295 mOsm/L adjusted with sucrose).
The intracellular solution comprised 120 mM KCl, 10 mM EGTA, 10 mM HEPES, 3 mM MgATP (pH 7.2 adjusted with KOH,
osmolarity 285 mOsm/L adjusted with sucrose).
Recordings were obtained at room temperature.

\paragraph*{Patching protocol}
On the day of recording, thick walled patch pipettes were pulled to a resistance of $\sim$ 3.5 M$\Omega$ and subsequently fire polished.
Whole cell patch clamp with seals of $>$ 5 G$\Omega$ were obtained in all cases, with series resistances $<$ 10 M$\Omega$, compensated for at the 70\% level.
Cells were then held at -60 mV for the duration of the application of the discrimination scheme with the control signal $u(t)$ applied to the analog input of the 470 nm wavelength channel of a CoolLED pe-800 via a National Instruments multifunction I/O module.
Responses to the control signal at each algorithm iteration were collected 5 times and averaged to reduce variability.
Consecutive applications of the control signal were spaced 2s apart whilst keeping the CoolLED off in between to prevent opsin desensitization.

\section{Evolution of parameters and control inputs}

\begin{figure}[h!]
    \centering
    \includegraphics[width=0.5\columnwidth]{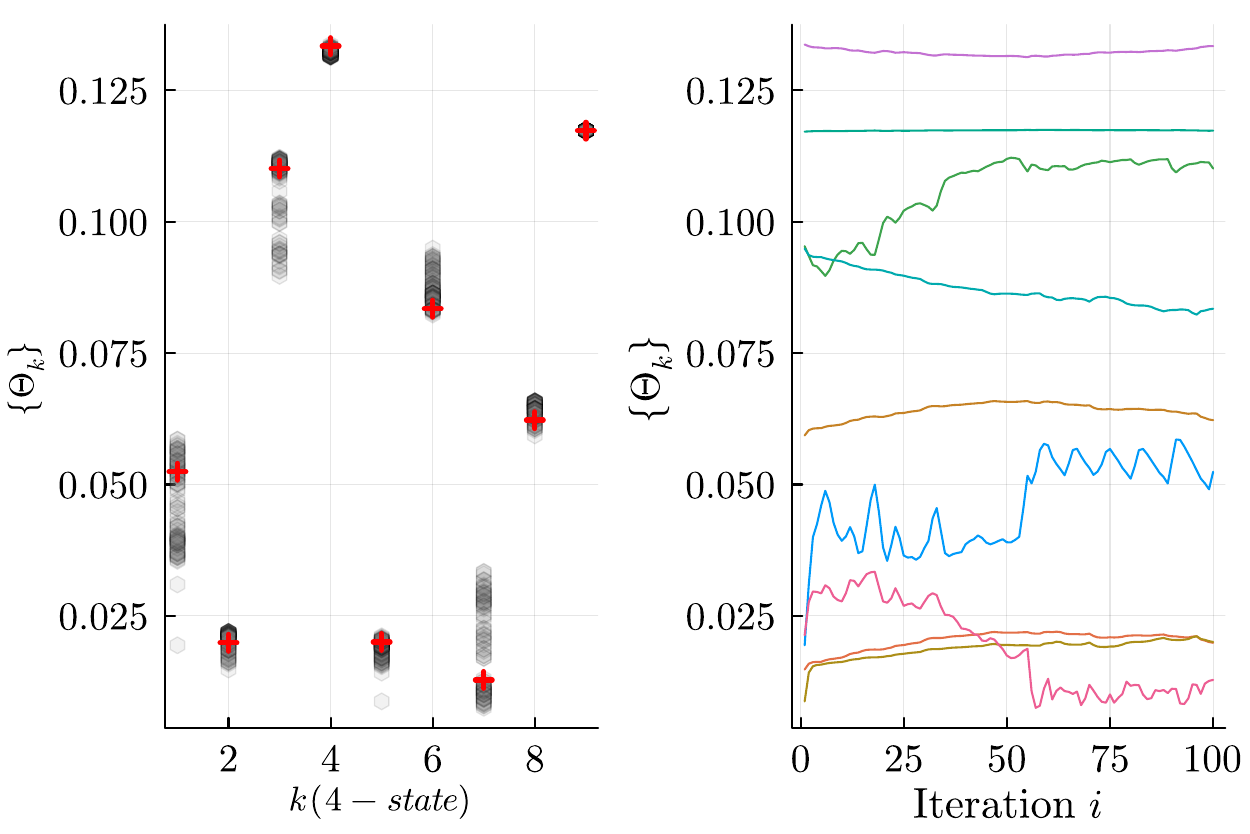}
    \caption{Evolution of the parameters for the 4-state model during the experiment shown in Fig.~\ref{fig:exp} (I). 
    Left panel: Parameters for all iterations (gray) and for the final iteration (red plus symbols). 
    Right panel: Trajectories of the parameters over the iterations. 
    Note that, for the sake of readability, we transformed
    $\theta_8 \rightarrow |\theta_8|$ and 
    $\theta_9 \rightarrow |\theta_9|/8$\,.}
    \label{fig:case_11}
\end{figure}

\begin{figure}[h!]
    \centering
    \includegraphics[width=0.5\columnwidth]{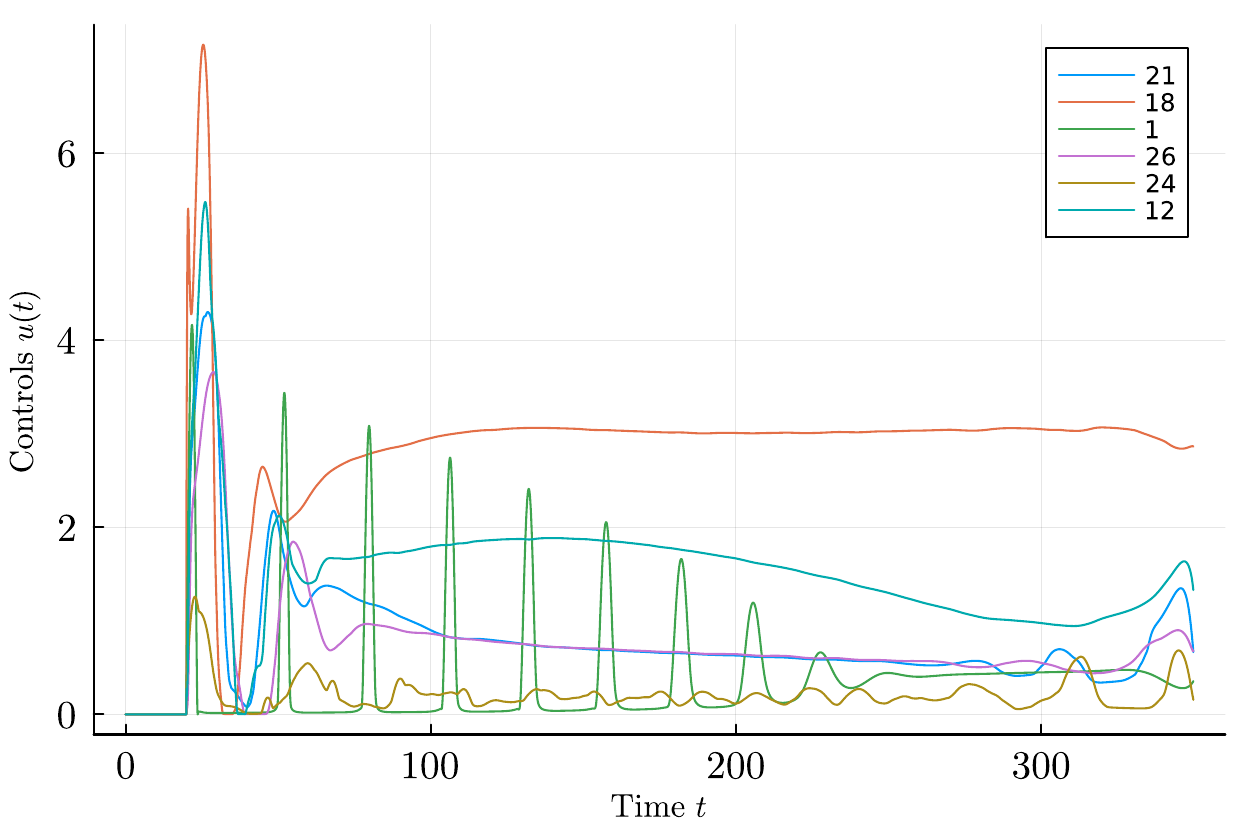}
    \caption{Comparison of the different control inputs obtained at each iteration of the discrimination process for the experiment shown in Fig.~\ref{fig:exp} (I). We used a $k$-means clustering algorithm with 6 groups. The solid lines give the centers of each cluster. The number of control inputs in each group is given in the legend. Note that many control inputs with different sequences of peaks are contained in the gold cluster. }
    \label{fig:case_11b}
\end{figure}

Figure \ref{fig:case_11} shows the evolution of the parameters of the 4-state model corresponding to Fig.~\ref{fig:exp} (I) with the number of iteration of the algorithm.
To show the diversity of control inputs, in Fig. \ref{fig:case_11b}, we performed a clustering of the signals generated at each iteration of the experiment in Fig.~\ref{fig:exp} (I).



\section{Identifiability of Opsin Models}\label{app:id}
The parameters of a model, given an input and an output, are said to be globally structurally identifiable if they can be uniquely recovered from the input-output relation~\cite{grewal1976identifiability}. 
Various algorithms have been proposed to determine the structural identifiability~\cite{bellu2007daisy,villaverde2016structural,dong2023differential}. 
Given the small number of prameters and variable in the 3-state model, one can assess the identifiability by hand. 
Following Ref.~\cite{meshkat2014finding}, using Eqs~(\ref{3s1})-(\ref{3s2})\, one obtains the input-output relation,
\begin{align}
\ddot{y} &= \theta_1\theta_2\theta_4 u(t)^2 - \theta_1\theta_2 y\, u(t)^2 - \theta_1\theta_3\theta_4 u(t)^3 + \theta_1\theta_3 y\, u(t)^3 + \theta_1\dot{y}\,u(t)^2 - \theta_2^2 y\, u(t) \\ \nonumber
&+ \theta_2\theta_3 y\, u(t)^2 - \theta_2 y\,\dot{u}(t) + \theta_3\dot{y}\,u(t)^2 - \dot{y}\,\dot{u}(t) + \ddot{y}\,u(t)\,.
\end{align}
The uniqueness of the parameters can then be shown by contradiction. 
Indeed, let us assume that there exists two sets of parameters $(\theta_1,\theta_2,\theta_3,\theta_4)$ and $(\phi_1,\phi_2,\phi_3,\phi_4)$ that produce the same output $y(t)$ given a diverse enough input $u(t)$\,.
Then taking the difference the input-output relation for the sets of parameters, one obtains the system of equations,
\begin{subequations}
\begin{align}
    \theta_1\theta_2\theta_4 - \phi_1\phi_2\phi_4 &= 0\,,\\
    \theta_1\theta_2 - \phi_1\phi_2 &= 0\,, \\
    \theta_1\theta_3\theta_4 - \phi_1\phi_3\phi_4 &= 0\,,\\
    \theta_1\theta_3 - \phi_1\phi_3 &= 0\,,\\
    \theta_1 - \phi_1 &= 0\,,\\
    \theta_2^2 - \phi^2 &=0 \,,\\
    \theta_2\theta_3 - \phi_2\phi_3 &= 0\,,\\
    \theta_2 - \phi_2 &= 0\,, \\
    \theta_3 - \phi_3 &= 0\,.
\end{align}
\end{subequations}
Solving the latter, one obtains $\theta_i = \phi_i$ for $i=1,2,3,4$\,, which contradicts the assumption.
A similar analysis can be performed for the 4- and 6-state models. However, it is a tedious procedure to do by hand. 
Instead, we checked the identifiability of both models using the algorithm presented in Ref.~\cite{dong2023differential}\,.

\section{Local Convergence under Strong Assumptions}\label{app:conv}
The goal of this is not to provide theoretical convergence guarantees, but rather to introduce a new method to discriminate models that outperform the standard procedures currently used in the literature, at least for electrophysiology and opsin modelling.

In most real world applications the hope of obtaining theory-based guarantees is vain (see the main text discussion on identifiability). Investigating the weakest set of assumptions guaranteeing convergence is left for future work. However, we would still like to show that the scheme is mathematically sound. In the following, we show the convergence of the our scheme under more stringent assumptions.

First, we assume full observability, i.e., the dimensions of the two models are the same and $\bm H_k$ is the identity map for $k=1,2$.
We assume that the ground truth (reference model) is an instance of model 1: 
\begin{equation}
\dot{\bm X}^* = \bm F_1(\bm X^*, \bm u, \bm \Theta^*) \,,
\label{eq:ref_model}
\end{equation}
with a fixed set of parameters $\bm \Theta^*$.

A model is locally identifiable under generic control if
\begin{align}
\mathcal{U}_{\,\rm fail} =\{\,\bm u :\;& \bm X(t,\bm \Theta, \bm u) =\bm X(t,\tilde{\bm \Theta}, \bm u)\,,\; \forall t \nonumber\\
&\text{ and }\bm \Theta \neq \bm \tilde{\bm \Theta} \,\}\subset \mathcal{U}\,,
\end{align} 
has measure zero. Or stated differently, $\forall \bm u \in \mathcal{U}_{\,\rm fail}^{\;c}$
\begin{equation}
\bm X(t,\bm \Theta, \bm u) =\bm X(t,\tilde{\bm \Theta}, \bm u) \Longrightarrow \bm \Theta =\bm \tilde{\bm \Theta}\,,
\end{equation} 
so that we are almost certain that trajectories are uniquely specified by the parameters once ${\bm u}$ is fixed.
We will first assume that the reference model is locally identifiable at $\bm \Theta^*$ under generic control, and later investigate the more general case.

We further assume that the second model is fixed, but captures the reference system rather accurately in the sense that
\begin{align}
\dot{\bm X}_2  &= \bm F_1(\bm X_2, \bm u, \bm \Theta^*) + \varepsilon\,\bm{\delta F}(\bm X_2, \bm u, \bm \Theta_2)\,,
\end{align}
where $\varepsilon=\sup_{\bm x,\bm u}\|\bm F_1(\bm x,\bm u)-\bm F_2(\bm x, \bm u)\|$ and $\bm{\delta F}(\bm X_2, \bm u, \bm \Theta_2)$ is the normalized deviation from the reference dynamics.

In order to facilitate the analysis, we simplify our scheme to
\begin{subequations}
\begin{align}
\bm u^{(k+1)}&={\rm argmax}_u\, \mathcal{J}(\bm \Theta_1^{(k)},\bm u)\,,\label{eq:u_step} \\
\bm \Theta_1^{(k+1)}&=\bm \Theta_1^{(k)}- \eta \nabla_\theta\,\mathcal{L}(\bm \Theta_1^{(k)}, \bm u^{(k+1)})\,,\label{eq:update}
\end{align}
\end{subequations}
where $\eta$ is a fixed learning rate. We label the optimal control obtained with $\bm \Theta^*$ by $\bm u^*$.  The simplifications are that: i) the mismatch maximization is now against a fixed target as is typically invoked by the T-optimality criterion~\cite{atkinsonDesignExperimentsDiscriminating,poncedeleon1991priorInfo}, and ii) the parameter update is not carried out via Adam, but with the simplest possible (one step) update rule. 
As in the main text, we will assume that there is no uncertainly on the initial state of the system.

As explained in the main text, we use the Euclidean distance to compare prospective models to the data (here reference model), thus
\begin{align}
\mathcal{L}\big(\bm \Theta^{(k)}, \bm u^{(k+1)}\big) = \frac{1}{2}\int_0^T& \big\|\bm X_1(t,\bm u^{(k+1)})
-\bm X^{*}(t,\bm u^{(k+1)})\big\|^2{\rm d}t\,.
\end{align}
Note that
\begin{align}
\mathcal{L}(\bm \Theta^*,\bm u)=0\,, \forall \bm u\,, \\
\nabla_\theta\mathcal{L}(\bm \Theta^*,\bm u^*)=0\,.
\end{align}

Expanding the gradient of the loss function $\mathcal{L}$ in Taylor series gives
\begin{align}
&\nabla_\theta\mathcal{L}(\bm \Theta_1,\bm u) =  \nabla_\theta\mathcal{L}(\bm \Theta^*,\bm u^{*})+ \nabla^2_{\theta\theta}\mathcal{L}(\bm \Theta^*,\bm u^*)(\bm \Theta_1-\bm \Theta^*)\nonumber\\
&+ \nabla^2_{\theta u}\mathcal{L}(\bm \Theta^*,\bm u^*)(\bm u-\bm u^*) + \text{(higher order terms)}\,.\label{eq:L_taylor}
\end{align}
The first and third terms on the right-hand side of Eq.~\eqref{eq:L_taylor} vanish, so we are left with 
\begin{equation}
\mathcal{H}:=\nabla_{\theta\theta}^2\mathcal{L}(\bm \Theta^*,\bm u)=\int_0^T \bm S^*(t)^{\top}\bm S^*(t)\,{\rm d}t\,,\label{eq:H}
\end{equation}

where  $\bm S^*(t):=\nabla_{\bm \theta} \bm X^*$ is the state sensitivity. It follows from differentiating \eqref{eq:ref_model} that
\begin{align}
\dot{\bm S}^* = \partial_x\bm F_1(\bm X^*,\bm u,\bm \Theta^*)\bm S^*
+  \partial_\theta\bm F_1(\bm X^*,\bm u,\bm \Theta^*)\,.\label{eq:S_dyn}
\end{align}

For self-consistency, let us state the following convergence theorem, which can be found in Ref.~\cite{ortega2000iterative}.

\textbf{Ostrowski Theorem:} Suppose that $\bm G: \mathcal{D} \subset \mathbb{R}^n \rightarrow \mathbb{R}^n$ has a fixed point $\bm x^* \in {\rm int}(\mathcal{D})$ and is Frechet-differentiable at $\bm x^*$. If the spectral radius of $\bm G'(\bm x^*)$ satisfies $\rho(\bm G'(\bm x^*)) < 1$, then $\bm x^*$ is a point of attraction
of the iteration 
\begin{equation}
\bm x_{k+1} = \bm G\, \bm x_k\,,\hspace{5pt} k=0,1,2,\ldots
\end{equation}
Combining Eqs.~\eqref{eq:update}, \eqref{eq:L_taylor} and \eqref{eq:H} one gets 
\begin{equation}
\bm \Theta_1^{(k+1)} = \big(\mathbb{I}- \eta \mathcal{H}\big)\, \bm \Theta_1^{(k)}\,.
\end{equation}
Thus, according to \textit{Ostrowski Theorem}, $\bm \Theta^*$ is a point of attraction if $\mathcal{H}$ is positive-definite, as it implies $\rho(\mathbb{I}- \eta \mathcal{H}) < 1$. We will now show that local identifiability of the reference model implies that $\mathcal{H}$ is positive-definite.

By construction, $\mathcal{H}$, defined in Eq.~\eqref{eq:H}, is a Gramian matrix. Therefore, it is at worst positive semi-definite. Now, suppose that $\mathcal{H}$ is singular and define $\delta {\bm X} = \bm X - \bm X^*$ and $\delta \bm \Theta = \bm \Theta - \bm \Theta^*$.
To first order, 
\begin{equation}
\int_0^T \delta \bm X^\top \delta \bm X\, {\rm d}t =  \delta \bm \Theta^\top \,\mathcal{H}\, \delta \bm \Theta\, .
\end{equation}

Thus, if $\mathcal{H}$ is singular, $\exists\, \bm v \in {\rm ker}(\mathcal{H})$ with $\|\bm v\| > 0$ so that 
\begin{equation}
\bm X(t,\bm \Theta^*, \bm u) =\bm X(t,\bm \Theta^* + \bm v, \bm u)\,, \forall t\,, \text{ with } \bm u\in\mathcal{U}_{\rm fail}^{\,c}\,,
\end{equation}
which contradicts local identifiability.
From the generic control assumption, the control obtained from Eq.~\eqref{eq:u_step} falls almost certainly in $\mathcal{U}_{\rm fail}^{\,c}$. Therefore, the scheme convergences locally, at least linearly \cite{ortega2000iterative}.

Now, what if we do not know whether or not the model in locally identifiable under generic control? In order to study this case,  let us simplify the problem further and assume that model 2, with fixed $\bm \Theta_2$, shares the same structure as model 1, so that
\begin{equation}
\dot{\bm X}_2 = \bm F_1(\bm X_2, \bm u, \bm \Theta_*+\bm \Theta_2)\,,
\end{equation}
where we assume that $\bm \Theta_{2,i}\neq0\,,\forall i$. Then the dynamics of $\bm\xi := \bm X_1 - \bm X_2$ becomes
\begin{align}
\dot{\bm \xi} &\approx \partial_x\bm F(\bm X^*,\bm u,\bm \Theta^*)\bm \xi +  \partial_\theta\bm F(\bm X^*,\bm u,\bm \Theta^*)\widetilde{\bm \Theta}\,,
\end{align}
with $\widetilde{\bm \Theta} := \bm \Theta_1 - \bm \Theta_2 -\bm \Theta^*$ and where we have dropped the subscript of $\bm F$. This dynamics is almost identical to the one for $\bm S^*$, given by Eq.~\eqref{eq:S_dyn}. Solving for $\bm S^*$ and $\bm \xi$ gives
\begin{align}
\bm S^*(t) &= \int_0^t\bm \Phi^*(\tau,0)\,\partial_\theta\bm F(\bm X^*,\bm u,\bm \Theta^*)\,{\rm d}\tau\,,\\
\bm \xi(t) &\approx \int_0^t\bm \Phi^*(\tau,0)\,\partial_\theta\bm F(\bm X^*,\bm u,\bm \Theta^*)\,\tilde{\bm \Theta} \,{\rm d}\tau\,,
\end{align}
where $\bm \Phi^*$ is state-transition matrix based on the reference trajectory and parameters. Thus, $\bm \xi(t) \approx \bm S^*(t) \,\widetilde{\bm\Theta}$.

To summarize, our control step \eqref{eq:u_step} tries to maximize
\begin{equation}
\mathcal{J}(\bm \Theta_1^{(k)},\bm u) \approx \widetilde{\bm\Theta}^\top\,\mathcal{H}(\bm u)\,\widetilde{\bm\Theta}\,,\label{eq:simpl_cont}
\end{equation}
and the question is: does this guarantee that $\mathcal{H}(\bm u)$ is positive-definite? We think that there are no ways to proof formally without deeper assumptions on the model's structure. In the simplified case shown in Eq.~\eqref{eq:simpl_cont}, it is easy to construct a counterexample by setting $\bm \Theta_1 = \bm \Theta_2 +\bm \Theta^*$ in which case all information is lost. Nonetheless, in most cases, it seems to promote the information about the model's parameters. In fact, we note that Eq.~\eqref{eq:simpl_cont} resembles to a Rayleigh quotient (without the normalizing factor), which is bounded by the smallest and largest eigenvalues of $\mathcal{H}$. Thus, maximizing the mismatch between the two models tends to increase the large eigenvalue by pushing it up from below. However, this procedure has no explicit effect on the smallest eigenvalue. So, in general, there is no guarantee that the mismatch-maximizing control will provide $\mathcal{H}$ positive-definite, but it is a sensible candidate for this as it increases the convergence rate.

During the in silico experiments that we did prior to this convergence analysis, we noticed that the scheme would sometimes get `stuck' in the vicinity of a control $\bm u_{\rm stuck}$. With the present analysis, we understand this as follows: the scheme gets into a place in the parameter space where the information it gets (almost) vanishes. We addressed this problem empirically by implementing the additional constraint that the current control should be as different as possible to previous control profiles. Furthermore, a gradient-descent scheme does not need information about every parameter at each iteration to converge. Coordinate descent, for instance, illustrates this point. Thus, as long as the scheme does not actively avoid a search region, it is reasonable to expect it to converge.

\end{document}